\documentclass[12pt]{article}
\usepackage{amsmath}
\usepackage{amssymb}
\usepackage{latexsym}

\usepackage{enumerate}

\newtheorem{theorem}{Theorem}[section]
\newtheorem{proposition}[theorem]{Proposition}

\newtheorem{lemma}[theorem]{Lemma}

\newtheorem{definition}[theorem]{Definition}
\newtheorem{notation}[theorem]{Notation}

\newenvironment{proof}{\noindent{\sc Proof.}}{\quad\qed\medskip}

\newcommand{\qed}{\quad\lower0.05cm\hbox{$\Box$}}

\newcommand{\Z}{{\mathbb Z}}
\newcommand{\K}{{\cal K}}

\newcommand{\Ks}{{\cal K}^*}

\newcommand{\R}{{\mathbb R}}
\newcommand{\G}{{\cal G}}
\newcommand{\Pa}{{\cal P}}
\newcommand{\adic}{^{\raise0.05cm\hbox{\scriptsize $\wedge$}}}

\newcommand{\wH}{\widehat H}
\newcommand{\wB}{\widehat B}
\newcommand{\wA}{\widehat A}

\newcommand{\Hom}{{\rm Hom}}
\newcommand{\Aut}{{\rm Aut}}
\newcommand{\Out}{{\rm Out}}

\newcommand{\ident}{{\rm id}}

\newcommand{\Inn}{{\rm Inn}\,}
\newcommand{\pinn}{{\rm pInn}\,}
\newcommand{\dom}{{\rm dom}\,}
\newcommand{\cod}{{\rm cod}\,}
\newcommand{\Dom}{{\rm Dom}\,}
\newcommand{\cf}{{\rm cf}\,}
\newcommand{\PP}{{\mathbb P}}

\newcommand{\GG}{{\mathbb G}}
\newcommand{\E}{{\mathbb E}}

\newcommand{\card}[1]{|#1|}
\newcommand{\norm}[1]{||#1||}
\newcommand{\Dl}{{\rm Dl}_\lambda}

\def\a{\alpha}
\def\b{\beta}
\def\g{\gamma}
\def\d{\delta}
\def\e{\varepsilon}
\def\l{\lambda}
\def\lp{\lambda^+}
\def\blp{\widetilde{\lambda^+}}
\def\bu{\widetilde{u}}

\def\k{\kappa}

\def\arr{\longrightarrow}
\def\iff{\Longleftrightarrow}

\def\Diam{\diamondsuit}

\def\restr{\upharpoonright}

\newcommand{\downarrowright}[1]{\downarrow
\rlap{\raise0.1cm\hbox{$\scriptstyle{#1}$}}}

\newcommand{\downarrowleft}[1]{\rlap{\kern-0.2cm
\raise0.1cm\hbox{$\scriptstyle{#1}$}}\downarrow}

\newcommand{\uparrowright}[1]{\uparrow
\rlap{\lower0.1cm\hbox{$\scriptstyle{#1}$}}}

\newcommand{\uparrowleft}[1]{\rlap{\kern-0.2cm
\lower0.1cm\hbox{$\scriptstyle{#1}$}}\uparrow}

\begin{document}
\setcounter{page}{1}
\title{\vspace*{-1cm} Large localizations of finite simple groups}

\author{
{\sc R\"{u}diger G\"{o}bel, Jos\'e L. Rodr\'{\i}guez and Saharon Shelah}
\thanks{The first and the third authors are supported by the project No. G
0545-173, 06/97
of the German-Israeli Foundation for Scientific Research \& Development.
The second author is supported by DGES grant
PB97-0202 and DAAD.
GbSh 701 in Shelah's list of publications.
}
}
\date{\today}
\maketitle


\begin{abstract}
A group homomorphism $\eta: H\to G$ is called a {\it localization}
of~$H$ if every homomorphism $\varphi : H\to G$ can be `extended
uniquely' to a homomorphism $\Phi :G\to G$ in the sense that $\Phi
\eta = \varphi$.

Libman showed that a localization of a finite group need not be
finite. This is exemplified by a well-known representation $A_n\to
SO_{n-1}(\R)$ of the alternating group $A_n$, which turns out to
be a localization for $n$ even and $n\geq 10$. Dror Farjoun asked
if there is any upper bound in cardinality for localizations of
$A_n$. In this paper we answer this question and prove, under the
generalized continuum hypothesis, that every non abelian finite
simple group $H$, has arbitrarily large localizations. This shows
that there is a proper class of distinct homotopy types which are
localizations of a given Eilenberg--Mac Lane space $K(H,1)$ for
any non abelian finite simple group~$H$.

\end{abstract}

\setcounter{section}{-1}


\section{Introduction}

\setcounter{equation}{0}

One of the current problems in localization of groups is to
decide what algebraic properties of~$H$ can be transferred to
$G$ by a localization $\eta: H\to G$.
Recall that $\eta: H\to G$ is a {\it localization}
if every homomorphism $\varphi: H\to G$ in the diagram
\begin{equation}
\label{localization}
\begin{array}{rrl}
H & \stackrel{\eta}{\longrightarrow}& G \\
\varphi\downarrow & \swarrow\hspace{-0.3cm} \raise -0.3cm\hbox{$\Phi$}\\
G &&
\end{array}
\end{equation}
can be extended to a {\it unique} homomorphism
$\Phi: G\to G$ such that $\Phi \eta =\varphi$.
This in other words says that $G\cong L_\eta H$
where $L_\eta$ is the localization functor with
respect to $\eta$; see e.g. \cite{Cas94}.

For example, the properties of being an abelian group or a
commutative ring with~1 are preserved. Casacuberta, Rodr\'{\i}guez and
Tai \cite{CRT97} have found consequences of these facts for homotopical
localizations of abelian Eilenberg--Mac Lane spaces (see also
Casacuberta~\cite{Cas94}, Dror Farjoun \cite{Dro96}).

A further step leads to nilpotent groups. Dwyer and Dror Farjoun
showed that every localization of a nilpotent group of class~2 is
nilpotent of class~2 or less. A proof is given by Libman in
\cite{Li98}.
However, it is unknown if nilpotent groups are preserved under
localizations in general.

Another interesting problem is to find an upper bound
for the cardinalities of the localizations of a fixed group $H$.
It is easy to see that, if $H$ is finite abelian, then every localization
$\eta: H\to G$ is an epimorphism, hence $\card{G} \leq \card{H}$.
More generally, Libman has shown in \cite{Li98} that if $H$ is torsion
abelian then $\card{G}\leq \card{H}^{\aleph_0}$.
However, if $H$ is not torsion, this can fail.
Indeed, the localizations of~$\Z$ are precisely the $E$-rings,
and we know by Dugas, Mader and Vinsonhaler \cite{DMV87} that there exist
$E$-rings of arbitrarily large cardinality. Str\"ungmann sharpened
this result in \cite{St98} for almost free $E$-rings.

The first example observing that a localization of a finite group
need not be finite is due to Libman \cite{Lib98}. He showed that
the alternating group $A_n$ has a $(n-1)$-dimensional irreducible
representation $\varphi: A_n \to SO_{n-1}(\R)$ which is a
localization for any even $n\geq 10$. In the proof he uses that
$O_{n-1}(\R)$ is complete, $SO_{n-1}(\R)$ is simple, and the fact
that all automorphisms of~$SO_{n-1}(\R)$ are conjugation by some
element in~$O_{n-1}(\R)$. This also motivates our Definition~1.
Thus, Dror Farjoun asked about the existence of an upper bound for
the cardinality of localizations of~$A_n$. We give an answer to
this question in Corollary~3, which is a direct consequence of our
Main Theorem and Proposition~2.

In fact our result also holds for many other groups which we will
call {\it suitable} groups, see Definition~1. However, we shall
assume the generalized continuum hypothesis. GCH will be needed to
apply a new combinatorial principle $\Dl$ due to Hart, Laflamme
and Shelah \cite{HLS93}. If $\l$ is a successor of a regular
cardinal, then under GCH this principle holds and is a general
tool which allows us to manufacture objects of size $\lp$ from
approximations of size less than $\l$.

The group theoretical techniques derive from standard
combinatorial group theory and can be found in the book by Lyndon
and Schupp \cite{LS77}. In particular, we will use small
cancellation theory as in \cite{She80}, where Shelah solved
Kurosh's problem about the existence of Jonsson groups.

\bigskip
\noindent {\bf Main Theorem} {\it (Assume {\rm ZFC + GCH}) Let
$\l$ be the successor of an uncountable, regular cardinal and let
$\lp$ be its successor cardinal. Then any suitable group $H$ whose
automorphism group $\Aut(H)$ has size $<\l$ is a subgroup of a
group $\G$ of size $\lp$ with the following properties:
\begin{itemize}
\item[{\rm (a)}]
$\G$ is generated by copies of $H$.
\item[{\rm (b)}]
Any monomorphism $\phi: H\to \G$ is the restriction of an
inner automorphism of $\G$.
\item[{\rm (c)}]
$H$ has trivial centralizer in $\G$, i.e. if $x\in \G$ commutes
with all elements of $H$, then $x=1$.
\item[{\rm (d)}]
$\G$ is complete, i.e. all automorphisms of~$\G$ are inner
automorphisms and $G$ has trivial centre $Z(\G)$.
\item[{\rm (e)}] $\G$ is co-hopfian, i.e. any monomorphism $\G \to \G$
is an automorphism.
\item[{\rm (f)}] $\G$ is simple.
\end{itemize}
}

\bigskip
\noindent {\bf Definition~1}
{\rm Let $H$ be any group with
trivial center and view $H \leq Aut(H)$ as inner automorphisms of
$H$. Then $H$ is called {\it suitable} if the following conditions
hold:
\begin{enumerate}
\item
$H$ is torsion, i.e. all elements have finite order.
\item If $H_1\leq Aut(H)$ and $H_1\cong H$ then $H_1=H$.
\item $H$ is complete in $\Aut(H)$, i.e. any automorphism of $H$
extends to an inner automorphism of $\Aut(H)$.
\end{enumerate}
}

We want to show that all  non abelian finite simple groups are
suitable. It only remains to show conditions 2 and 3 in the above
definition. As a consequence of the classification of finite
simple groups the Schreier conjecture holds for all finite
groups~$H$, hence the outer automorphism group $\Out(H)$ is always
soluble. If $H$ is also simple and non abelian then $H$ is a
characteristic subgroup of $\Aut(H)$. Hence any copy of $H$ in
$\Aut(H)$ is $H$ and any automorphism of $H$ is induced by an
automorphism of $\Aut(H)$. However, $\Aut(H)$ is complete by
Burnside, see \cite[p.~399]{Rob}. This is the first part of the
following:

\medskip
\noindent
{\bf Proposition~2}
{\it
{\rm (a)} All non abelian finite simple groups are suitable.

\noindent {\rm (b)} The non abelian, finite, simple and complete
groups are precisely the following sporadic groups: $M_{11},
M_{23}, M_{24}, Co_{3}, Co_{2}, Co_{1}, Fi_{23}, Th, B, M, J_1,
Ly, Ru, J_4.$ }

\medskip
\noindent
Part (b) follows by inspection of the list of finite simple groups. It is
interesting to know when $H=\Aut(H)$ because in this case
our proof becomes visibly simpler.

We finally obtain an answer to Dror Farjoun's question.

\medskip
\noindent
{\bf Corollary~3}
{\it (Assume {\rm ZFC + GCH}.) Any suitable simple group
has localizations of arbitrarily large cardinality. }

\medskip

The localization $A_n\to SO_{n-1}(\R)$ for any even $n \geq 10$
induces a map between Eilenberg--Mac Lane spaces $K(A_n, 1)\to
K(SO_{n-1}(\R),1)$ which turns out to be a localization in the
homotopy category; \cite{Lib98}. This is the first example of a
space with a finite fundamental group which admits localizations
with an infinite fundamental group. Corollary~3 yields then the
following extension:

\medskip
\noindent
{\bf Corollary~4}
{\it (Assume {\rm ZFC + GCH}.) Let $H$ be a suitable
simple group. Then $K(H,1)$ has localizations with arbitrarily
large fundamental group.
}

\medskip
\noindent Other constructions in homotopy theory based on
large-cardinal principles were described by Casacuberta, Scevenels
and Smith in \cite{CSS98}.

\medskip
\noindent

\bigskip\noindent
{\bf Acknowledgements.}  We would like to thank Emmanuel Dror Farjoun for
explaining the above problem to us and pointing out further interesting
connections between group theory and homotopy theory.


\section{Proof of corollaries}

\setcounter{equation}{0}

It is known that the generalized continuum hypothesis ensures the
existence of arbitrarily large cardinals $\lambda$ which satisfy
$\Dl$. For instance, successor cardinals $\l=\kappa^+$ of regular
cardinals $\kappa>\aleph_0$ satisfy $\Dl$. Hence, in order to
prove Corollary~3, we next show that the inclusion $\eta: H
\hookrightarrow \G$ given in the Main Theorem is a localization.
Suppose that $\varphi: H\to \G$ is a homomorphism. We have to show
that there is a unique homomorphism $\Phi:\G\to \G$ such that
$\Phi \eta =\varphi$. If $\varphi = 0$ then $\Phi = 0$ makes the
diagram (\ref{localization}) commutative. To see that it is
unique, we note that $H$ is in the kernel  $K$ of $\Phi$ and $K =
\G$ by simplicity, hence $\Phi = 0$. Now suppose that $\varphi\neq
0$. Since $H$ is simple we have that $\varphi(H)\leq \G$ and
$H\cong \varphi(H)$. Thus by part (b) of the Main Theorem there is
an element $y\in \G$ such that $\varphi= y^* \restr H$ where $y^*$
is the inner automorphism of~$\G$ given by $y^*(x)= y^{-1}x y$ for
all $x \in \G$, and $y^*\restr H$ denotes its restriction on $H$.
Hence $\Phi=y^*$ satisfies $\Phi \eta =\varphi$. Suppose that
$\Phi': \G\to \G$ is another homomorphism such that $\Phi' \eta
=\varphi$. Then $\Phi'\neq 0$. Now since $\G$ is simple by (f) and
co-hopfian by (e) we have that $\Phi'$ is an automorphism. Thus,
by completeness (d) we have $\Phi'=g^*$ for some $g\in \G$. Now
$(g^{-1} y)^* :\G\to \G$ fixes all elements of $H$. Hence, by part
(c), we finally obtain that $g^{-1} y=1$ and thus $\Phi=\Phi'$ as
desired.

Recall from \cite{Dro96} or \cite{CRT97} that a map $f: X\to Y$
between two connected spaces is a {\it homotopical localization}
if $Y$ is $f$-local, i.e. if the map of pointed function spaces $$
{\rm map_*\,} (Y,Y) \to {\rm map_*\,} (X,Y) $$ induced by
composition by $f$ is a weak homotopy equivalence. As in the case
of groups this says that $Y\simeq L_fX$, where $L_f$ is the
localization functor with respect to $f$.

It turns out that the homotopical localizations of the circle
$S^1=K(\Z,1)$ are precisely Eilenberg--Mac Lane spaces $K(A,1)$
where $A$ ranges over the class of all $E$-rings\linebreak
\cite[Theorem~5.11]{CRT97}. And therefore this class is proper
(i.e. it is not a set) in view of the result in \cite{DMV87}.
Recall that an $E$-ring is a commutative ring $A$ with identity
which is canonically isomorphic to its own ring of additive
endomorphisms. Corollary~4 claims that a similar situation holds
for $K(H,1)$ if $H$ is a suitable group. However, in this case,
other localizations of $K(H,1)$ which are not of the form $K(G,1)$
may exist.

Corollary~4 follows from the fact
that every localization of groups $H\to \G$ gives rise to
a localization of spaces $K(H,1)\to K(\G,1)$.
This holds because, for arbitrary groups $A$ and $B$,
the space ${\rm map_*\,} (K(A,1), K(B,1))$
is homotopically discrete and equivalent to the set $\Hom(A,B)$.

 We recommend the reader not so familiar with the basic facts on
homotopy theory to consult, for instance, \cite{Whi78} or
the monograph by Aubry \cite{A} (lectures of a DMV seminar
by Baues, Halperin and Lemaire) as well as Dror Farjoun's
exposition \cite{Dro96} for homotopical localizations.


\section{ Sketch of the proof of Main Theorem}
\label{Sketch}
\setcounter{equation}{0}

The group $\G$ in our main theorem has cardinality $\lp$, where
$\l^+$ is the successor cardinal of a regular cardinal~$\l$. It
will be the union of a $\lp$-filtration $\{\G_\a : \a<\lp\}$, that
is, $$ \G= \bigcup_{\a<\lp} \G_\a $$ where $\G_0=Aut(H)$ and $|\G_\a|=\l$
for $0 \neq \a< \lp$. Moreover $|\G_{\a+1}\setminus \G_\a|=\l$ as
well. More precisely, at each step $\a$, $\G_\a$ will be the union
of a certain directed system $\GG_\a$ of size $|\GG_\a|=\l$, where
each subgroup $G\in \GG_\a$ has cardinality $< \lambda$. It is
convenient to give a linear order to all elements which are added
in the construction of $\G$, so we view any $G \in \GG_\a$ as a
subset of $$ \blp := \{\l\a' +i : \a'\in \lp,\, i<\l\}. $$ This
set equals $\lp \times \l$ with the lexicographical order on the
pair of ordinals, hence it is well-ordered and has cardinality
$|\blp| = \lp$. Hence $\blp$ will be the universe for~$\G$. All
groups $G\in \GG_\a$ ($0\neq\a<\lp$) --which can be thought as
``approximations'' to $\G$-- will satisfy properties (a), (b) and
(c) of the Main Theorem, and moreover $G\subseteq \l\zeta_\a$ for
some unbounded, strictly increasing sequence of ordinals
$\zeta_\a< \lp$. The construction of these approximations is given
in Section~\ref{a-c} (only as groups). We show in
Section~\ref{upo} how they can be embedded into $\blp$; see
Proposition~\ref{new-Ks-dense-K}. The group structure of each $G$
is coded, up to a certain strong isomorphism, by a pair $p=(\g,u)$
with $\g\in \l$ and $u \subset \lp$ of cardinality $\card{u}<\l$
($u$ is the image of $G$ under the natural projection $\blp\to
\lp$). All such pairs $p$ are collected in a partially ordered set
$\PP$ which satisfies the requirements of the $\PP$-game defined
by Hart, Laflamme and Shelah in \cite{HLS93}, which we have
included in Appendix~A for convenience. A ``local'' construction
taking care of simplicity (f) is done in Proposition~\ref{Brit}.
The properties to be complete (d) and co-hopfian  (e) are treated
together in Section~\ref{Section-small-cancellation} and again in
Section~\ref{Section-small-cancellation-blp} inside the universe
$\lp$. Finally, in Section~\ref{PROOF} we explain how the two
players of the $\PP$-game carry out the ``global'' construction of
$\G$: player~I defines adequate ``density systems'' which make
some restrictions to player~II in his choice of the sequence
$\zeta_\a$ and the directed systems $\GG_{\a}$. These restrictions
will guarantee the desired properties of $\G$. The combinatorial
principle $\Dl$, which is a weak version of diamond
$\Diamond(\l)$, is used by player~II to ``predict'' in some sense
the moves of player~I.


\section{Some preliminaries on set theory, free products
with amalgamation and HNN extensions}
\label{set-theory}

\setcounter{equation}{0}

We recall some basic notions of set theory from \cite{Jech}.

The cofinality $\cf(\l)$ of an infinite cardinal $\l$ is the
minimal cardinal $\k$ for which there exists an increasing
sequence of ordinals $\l_i<\l$ of length $\k$
such that $\sup\{\l_i: i<\kappa\}=\l$.
Clearly, $\cf(\l)\leq \l$. A cardinal $\l$ is called {\it regular}
if $\cf(\l)=\l$.

Let $\l$ be an uncountable regular cardinal.
A subset $C\subseteq \l$ is a {\it club} if it is closed
under limits of length $<\l$ and if it is unbounded.
A subset $S\subseteq \l$ is {\it stationary} if $S\cap C\neq \emptyset$
for all clubs $C$ in $\l$. In particular every club is stationary.
Let $\lp$ be the successor cardinal of $\l$. We will use
that
$
C=\{\d <\lp:  \cf(\d)=\l \}
$
is a club of $\lp$. Moreover, to prove that the group
$\G$ is complete (Section~9.4) we need Fodor's Lemma,
see Jech \cite[p.~59, Theorem~22]{Jech}.
\begin{lemma}
\label{Fodor}
Suppose $\eta: S'\to \lp$ is a regressive function where $S'$ is a
stationary subset of $\lp$. Then there is a stationary $S \subseteq S'$
such that $\eta$ is constant on $S$.
\end{lemma}

Recall that the generalized continuum hypothesis GCH says that
for any cardinal $\k$ the successor $\k^+$ coincides with
$\card{\Pa(\k)} =2^\k$, where $\Pa(\k)$ denotes the set of subsets of~$\k$.
Moreover GCH implies that $\card{\Pa_{<\l}(\l)}=\l^{<\l}$ is $\l$,
where $\Pa_{<\l}(\l)$ denotes the set of subsets of $\l$ of
cardinality $<\l$.

In the next section we shall need two lemmas describing torsion
subgroups of free products with amalgamations and HNN extensions
(compare with Theorem~2.7, p.~187, Theorem~6.6, p.~212 and
Theorem~2.4, p.~185 in \cite{LS77}).

\begin{lemma}
\label{torsion-theorem}
Let $G^*=G *_{G_0}G'$ be the free product of $G$ and $G'$
amalgamating a common subgroup $G_0$. If $H'$ is a torsion
subgroup of $G^*$ then there exists a $y\in G^*$ such
that $(H')^y$ is contained in $G$ or in $G'$.
\end{lemma}
\begin{proof}
Suppose $H'$ is not contained in $G_0$, otherwise we are done. Fix
$h\in H'\setminus G_0$ and let $h$ be written in a canonical form
$$ h= g_1g_1' \cdots g_ng_n',$$ i.e. $g_i\in G\setminus G_0$ and
$g_i'\in G'\setminus G_0$ (with the convention that $g_1$ or
$g_n'$ may be trivial.) Note that this expression is not unique in
general, but its length $|h|$ only depends on $h$.

Suppose first that all elements of $H'$ have length 1.
If $H'$ contains two elements $g\in G\setminus G_0$
and $g'\in G'\setminus G_0$ then
$gg'\in H$ is a word written in a reduced normal form,
and has length 2, hence it is non torsion
by Theorem~2.6, p.~187 in \cite{LS77}.
This contradicts that $H$ is torsion.
Thus, all $g\in H'$ must belong
to the same factor, and we are done.

We can suppose that the first factor $g_1$ of $h$ is non trivial
and that $|h|$ is odd and $\geq 3$.
Choose $h$ such that $|h|$ is the minimum over all
lengths of elements in $H'$ of length $\geq 3$.
If $g_n'\neq1$ then we obtain a contradiction as before.
Thus $g_n'$ must be trivial and similarly $g_ng_1\in G_0$.
Therefore, $h$ has a normal decomposition of
the form $$h= g_1 g_1' \cdots \overline g_{n-1}'g_1^{-1},$$
for some $\overline g'_{n-1}\in G'\setminus G_0$ which
we rename by $g'_{n-1}$ to simplify the notation.

Similarly, all $h_1\in H'$ must have odd length. Moreover, if
$|h_1|\geq 3$ then $h_1$ has a normal form starting with $g_1$,
the first factor of $h$ and ending with $g_1^{-1}$. Of course, if
$|h_1|$=1 then $h_1\in G$.

Consider then $g_1^{-1}H' g_1$. If $h_1\in H'$ has length $\geq
3$, then $g_1^{-1}h_1g_1$ has length $|h|-2$, and its normal form
begins and ends with elements in $G'$. On the other hand, if
$h_1\in H'$ has length 1, then $g_1^{-1} h_1 g_1 \in G$. But,
arguing as before, we see that such elements of length 1 must be
actually in $G_0$.

Repeating this process we obtain a $y\in G^*$ such that $h= y g
y^{-1}$ where $g\in G\setminus G_0$ or $G'\setminus G_0$ depending
on $n$, if $n$ is odd or even, respectively. If any of the
elements in $(H')^y$ has length $\geq 3$, then we can conclude
that $g\in G_0$ which is a contradiction. Thus all elements in
$(H')^y$ must have length 1, and the first case applies, hence
$(H')^y$ is contained in some factor. This finishes the proof.
\end{proof}

\begin{lemma}
\label{torsion-theorem-HNN}
Let $G$ be any group, and $\phi: G_0\to G_1$ be an isomorphism
between two subgroups of $G$. Consider the HNN extension
$G^* = \langle G, t \,: \, t^{-1} h t= \phi(h), h\in G_0\rangle$.
If $H'$ is a torsion subgroup of $G^*$, then there exists a
$y\in G^*$ such that $(H')^y$ is contained in~$G$.
\end{lemma}
\begin{proof}
The proof is similar to that of Lemma~\ref{torsion-theorem}.
\end{proof}

The following lemma describes centralizers of torsion subgroups in
free products with amalgamation. An analogous lemma for HNN
extensions also holds.

\begin{lemma}
\label{nontrivial-centralizer} Let $G^*=G *_{G_0}G'$ be the free
product of $G$ and $G'$ amalgamating a common subgroup $G_0$. Let
$H'\leq G$ be a torsion subgroup and let $1\neq x\in {G^*}$ an
element which commutes with all elements of $H'$. Then either
$x\in G$ or $(H')^g \leq G_0$ for some $g\in G^*$.
\end{lemma}
\begin{proof}
Let $H'\leq G$ be a torsion subgroup, and
$x\in G^*$ such that $[x,H']=1$. If $x\in G$ we are done.
Suppose then that $x\not\in G$.
Express $x$ in a reduced normal form
$$x= g_1g_1' \cdots g_ng_n',$$
that is, $g_i\in G\setminus G_0$ and $g'_i\in G'\setminus G_0$,
where $g'_1\neq 1$ ($g_1$ or $g_n'$ may be trivial).
The relation $h^{-1} x ^{-1} h x =1$ yields the following
$$
h^{-1} {g'}_n^{-1} g_n^{-1} \cdots {g'}_1^{-1}
(g_1^{-1} h g_1) g'_1 \cdots g_n g'_n =1.
$$
By the normal form theorem for free products with amalgamation
\cite[Theorem~2.6 p~187]{LS77}, this is only possible if
$g_1^{-1} h g_1 \in G_0$ for
all $h\in H'$. This concludes the proof.
\end{proof}


\section{Approximations to $\G$}
\label{a-c}

\setcounter{equation}{0} In this section we fix a suitable group
$H$, we denote $\widehat H = \Aut(H)$ for shortness and view
$H\leq \widehat H$ as subgroup. We also fix an uncountable regular
cardinal $\l$ such that $|\wH|< \l$. As usual $C_{G'}(G)$ and
$N_{G'}(G)$ denote, respectively, the centralizer and the
normalizer of a subgroup $G$ in a group $G'$.

We next consider two particular families $\Ks \subseteq \K$ of
groups which will be used to approximate our group $\G$ of the
Main Theorem in the sense explained in Section~\ref{Sketch}.

First we prove a series of lemmas which culminate in
Proposition~\ref{Ks-dense-K}.
This states that the class of groups $\Ks$ is {\it dense} in~$\K$,
that is, given an infinite group $G\in \K$ there is a group $G'\in \Ks$ of
the same cardinality such that $G\leq G'$.
In particular, we provide an explicit algorithm to construct
groups satisfying conditions
(a), (b) and (c) of the Main Theorem.

We also show Proposition~\ref{Brit} which will be used to prove that
our final group $\G$ is simple.

\begin{definition}
\label{K-Ks}
{\rm
\begin{enumerate}
\item
$\K$ consists of all groups $G$ with $\card{G}<\l$
such that $H\leq \wH \leq G$, and any isomorphic copy of $H$ in $G$
has trivial centralizer in $G$. That is,
$$
\K= \{G: \mbox{ $\wH \leq G$, $\card{G}<\l$,
if $H\cong H'\leq G$, $x\in G$ with $[H', x]=1$,  then $x=1$}\},
$$
where $[x, y]= x^{-1}y^{-1}xy$.
\item
$\Ks$ is the class of all groups $G$ in~$\K$
which are generated by copies of~$H$ and
such that any isomorphism between $H$
and a subgroup of $G$ is induced by conjugation with an
element in~$G$.
\end{enumerate}
}
\end{definition}
Note that condition (1) in Definition \ref{K-Ks} is stronger than
property (c) of the Main Theorem. If we denote by $H(G)$ the {\it
$H$-socle} of~$G$, i.e. the subgroup of~$G$ generated by all
images of homomorphisms from $H$ into $G$, then $H(G)=G$ for all
$G\in\Ks$. We start observing that $\K$ is not empty.
This is obvious since $H$ is suitable and $|\wH|<\l$.
\begin{lemma}
\label{HinK}    $\wH$ is in  $\K$. \qed
\end{lemma}

\begin{lemma}
\label{free-groups-K}
If $G$ and $G'$ are in  $\K$ then $G*G' \in \K$.
\end{lemma}
\begin{proof}
Suppose that $H'\leq G* G'$ with $H'\cong H$ and $x\in G*G'$ such
that $[H',x]=1$.
By Lemma~\ref{torsion-theorem}, we can suppose that $(H')^y\leq
G$. Hence $[H',x]=1$ implies $[(H')^y , x^y]= 1$. In other words,
\begin{equation}
\label{commutator}
(h^y)^{-1} (x^y)^{-1}h^y x^y =1
\end{equation}
 for all $h\in H'$.
If $x^y \not\in G$, then express $(x^y)^{-1}$ in reduced form and
replace a first choice $h \neq 1$ by a different one if the first
$G$-factor of $(x^y)^{-1}$ in (\ref{commutator}) is cancelled by
$(h^y)^{-1}$. Hence the normal form of the left hand in
(\ref{commutator}) has non-trivial factors and the normal form
theorem for free products \cite[p.~175]{LS77}  gives a
contradiction, hence $x^y\in G$. But $G\in \K$, hence $x^y=1$ and
thus $x=1$, as desired. Note that the lemma also follows from
Lemma~\ref{nontrivial-centralizer}.
\end{proof}

\medskip
Lemma~\ref{HNN} below will take care of part (b) of the Main
Theorem. If $G$ is any group in $\K$ and $\phi:A\to B$ is an
isomorphism between two subgroups of $G$ isomorphic to $H$, we
want that $\phi$ is the restriction of an inner automorphism in
some extension $G\leq G_2\in \K$. This can be established using
HNN extensions as we next show. We will get this after several
steps.

\begin{lemma}
\label{HNN-1}
Let $G\in \K$ and $B\leq G$ be a subgroup isomorphic to $H$.
Then there is $G\leq G_1\in \K$ such that $\Aut(B) \leq G_1$.
\end{lemma}
\begin{proof}
Let $\wB=\Aut(B)$ and $N=N_G(B)$. If $\wB\leq G$ then $G_1=G$.
Suppose that $\wB\not \leq G$. Note that $N=G\cap \wB$ so we can
consider the free product with amalgamation $G_1:= G*_{N} \wB$. We
shall show that $G_1\in \K$. Let $H'\leq G_1$ be a subgroup
isomorphic to $H$ and $1\neq x\in G_1$ such that $[H',x]=1$. By
Lemma~\ref{torsion-theorem} we can suppose that $H'\leq G$ or
$H'\leq \wB$. Suppose that $H'\leq G$, the other case is easier.
Let $x=g_1g_2 \cdots g_n$ be written in a reduced normal form. If
$x=g_1 \in G$ then $x=1$ since $G\in \K$, and this is a
contradiction. Hence $x=g_1\in \wB\setminus N$. As in
Lemma~\ref{nontrivial-centralizer} we deduce that
$H'=(H')^{g_1}\leq N$, thus $H'=B$ since $B$ is suitable
(condition~2). Hence $g_1\in N$ and we get a contradiction. If $n=
2$, then we obtain $(H')^{g_1}= (H')^{g_2}= B =H'$. So both $g_1$
and $g_2$ are in $N$, which is a contradiction. Similarly, if
$n\geq 3$ we have that $g_2$ and $g_3$ are in $N$. This is again
impossible. This concludes the proof.
\end{proof}

By the previous lemma we can suppose that if $B \leq G\in \K$, and
if $B\cong H$, then $\wB \leq G$ as well.

Recall that two subgroups $A$ and $B$ of a group $G$ are conjugate
in $G$ if there is $g\in G$ such that $A=B^g$. In particular, if
$A, B \leq G$, $A\cong B\cong H$ and $\wA$, $\wB$ are conjugate in
$G$ then $A$ and $B$ are also conjugate. Indeed, if $g\in G$ such
that $g^*: \wA\cong \wB$, then the image $A^g\leq \wB$ is a
subgroup isomorphic to $B$, hence $A^g=B$, by condition~2 of
Definition~1.

\begin{lemma}
\label{HNN-2} Let $G\in \K$ and $B\leq \wB \leq G$, where $B\cong
H$. Suppose that $H$ and $B$ are not conjugate in $G$. Let $\phi:
\wH \to \wB$ be any isomorphism. Then, the HNN extension $$ G_1=
\langle G, t \,: \, t^{-1} h t= \phi(h) \mbox{\,\, for all $h\in
\wH$}\rangle $$ is also in $\K$.
\end{lemma}
\begin{proof}
Let $H'\leq G_1$ be a subgroup isomorphic to $H'$
and $1\neq x\in G_1$ such that $[H',x]=1$. By Lemma~\ref{torsion-theorem-HNN}
we can suppose that $H'\leq G$ already.
Let
$$
x=g_0 t^{\e_1}g_1t^{\e_2}  \cdots g_{n-1}t^{\e_n}g_n
$$
be written in a reduced form in~$G_1$, where $g_i\in G$ and there is no
subwords $t^{-1}g_it$ with $g_i\in \wH$ or $tg_it^{-1}$
with $g_i\in \wB $ (see \cite[p.~181]{LS77}).

If $x=g_0\in G$ then $x=1$, since $G\in \K$. This yields a contradiction.
Thus $n\geq 1$.
We have $[h,x]=h^{-1}x^{-1}hx=1$ for every $h\in H'$.
In other words, for a fixed $h\neq 1$, the following holds
\begin{equation}
\label{hxhx}
h^{-1}g_n^{-1} t^{-\e_n} \cdots t^{-\e_1}
(g_0^{-1}hg_0)
t^{\e_1} \cdots t^{\e_n}g_n =1.
\end{equation}
By the normal form theorem for HNN extensions
(\cite[p.~182]{LS77}), either $\e_1=1$ and $g_0^{-1} h g_0\in
\wH$, or $\e_1=-1$ and $g_0^{-1} h g_0\in \wB$. Suppose that
$\e_1=1$, the other case is analogous. Then $(H')^{g_0}\leq \wH$,
thus $(H')^{g_0}=H$, and we can replace in  (\ref{hxhx}) the
subword $t^{-1} (g_0^{-1}hg_0) t$ by $\phi(g_0^{-1}hg_0) \in B$.
Repeating the same argument we obtain that $\varepsilon_i=1$ or
$-1$. Hence one of the two possibilities holds depending on
$\varepsilon_2=1$ or $\varepsilon_2=-1$. So we have either $$
\ident : H_1 \stackrel{g_0^*}{\longrightarrow} \wH
\stackrel{\phi}{\longrightarrow} \wB
\stackrel{g_1^*}{\longrightarrow} \wH
\stackrel{\phi}{\longrightarrow} \cdots
\stackrel{g_n^*}{\longrightarrow} H_1 $$ or $$ \ident : H_1
\stackrel{g_0^*}{\longrightarrow} \wH
\stackrel{\phi}{\longrightarrow} \wB
\stackrel{g_1^*}{\longrightarrow} \wB
\stackrel{\phi^{-1}}{\longrightarrow} \cdots
\stackrel{g_n^*}{\longrightarrow} H_1. $$ In the first case we
have an isomorphism $g_1^* \phi: \wH
\stackrel{\phi}{\longrightarrow} \wB
\stackrel{g_1^*}{\longrightarrow} \wH$. Since $\wH$ is complete
there is $g\in \wH$ such that $g_1^* \phi=g^*$. This yields
$\phi=(g_1^{-1}g)^*$, i.e. $\wH$ and $\wB$ are conjugate, and thus $H$
and $B$ are conjugate, but this is impossible by hypothesis. In
the second case we have $g_1\in \wB$ by completeness. But, on the
other hand $g_1\not\in \wB$ since $x$ is written in a reduced
form. We conclude that $G_1$ is in $\K$, as desired.
\end{proof}

Let $\pinn(G)$ denote the set of {\it partial inner automorphisms},
which are the isomorphisms $\phi: G_1\to G_2$
where $G_1, G_2\leq G$ such that $\phi$ can be extended to an inner
automorphism of~$G$.
Hence  $\pinn(G)$ are all restrictions of conjugations to subgroups of~$G$.

\begin{lemma}
\label{HNN}
Let $G\in \K$. Let $A\leq \wA\leq G$ and $B\leq \wB\leq G$.
Let $\phi: A\cong B$ be any isomorphism.
Then, there is $G_2\in \K$ such that $\phi\in \pinn(G_2)$.
Moreover, $G_2$ can be obtained from $G$ by at most two successive HNN extensions.
\end{lemma}
\begin{proof}
If $\phi\in \pinn(G)$ we take $G_2=G$. Suppose that $\phi\not\in
\pinn(G)$. If $\wH$ and $\wA$ are conjugate, we take $G_1=G$.
Otherwise, we consider the HNN extension $$G_1=\langle G, t_1:
t_1^{-1} h t_1 = \phi_1(h)\mbox{\,\, for all $h\in \wH$}\rangle$$
where $\phi_1$ is any isomorphism between $\wH$ and $\wA$. By
Lemma~\ref{HNN-2} we know that $G_1\in \K$. Now, if $H$ and $B$
are conjugate in $G_1$, we take $G_2=G_1$. It follows
automatically that $\phi\in \pinn(G_1)$, since $\wH$ is complete.
If $\wH$ and $\wB$ are not conjugate in $G_1$, we consider a new
HNN extension $$G_2=\langle G_1, t_2: t_2^{-1} h t_2 =
\phi_2(h)\mbox{\,\, for all $h\in \wH$}\rangle$$ where $\phi_2$ is
any isomorphism between $\wH$ and $\wB$. Again $G_2 \in \K$ by
Lemma \ref{HNN-2}. In that case we have an isomorphism
$(t_2^{-1})^* \phi t_1{^*}: \wH\to \wH$, which equals $g^*$ for
some $g\in\wH$ by completeness. Thus $\phi= (t_2 g
t_1^{-1})^*\restr A$. This shows that $\phi \in \pinn(G_2)$, as
desired.
\end{proof}

The following two lemmas will be used to prove Lemma~\ref{HG=G}.
This will take care of part (a) of the Main Theorem.
Let $o(g)$ denote the order of $g\in G$.

\begin{lemma}
\label{ogog}
Let $G\in \K$ and suppose that $G'\in \K$ or $G'$ does not contain
any subgroup isomorphic to $H$.
Let $g\in G$ and $g'\in G'$ with $o(g)=o(g')$.
Then $(G*G')/N\in \K$
where $N$ is the normal subgroup of $G*G'$ generated by $g^{-1}g' \in G*G'$.
\end{lemma}
\begin{proof}
Let $\overline G=(G*G')/N$ and view $G$ and $G'$ as subgroups of
$\overline G$, respectively. Suppose that we have a subgroup
$H'\leq \overline G$ isomorphic to $H$ and $x\in \overline G$
such that $[H',x]=1$. By Lemma~\ref{torsion-theorem}
a conjugate of $H'$ is contained in~$G$ or $G'$. Hence
by hypothesis we can assume that $H'$ is already contained in~$G$.
Suppose that $x\neq 1$. By Lemma~\ref{nontrivial-centralizer}
it follows that either $x\in G$ or a conjugate of $H'$ is contained in
$\langle g \rangle$. In the first case $x=1$ since $G\in \K$, so contradiction,
and the second case is obviously impossible.
Thus $\overline G\in \K$ as desired.
\end{proof}

\begin{lemma}
\label{og}
If $g\in G\in \K$, then there is a group $\overline G\in \K$,
such that $G\leq \overline G$, with
$\card{\overline G}=\card{G} \cdot \aleph_0$ and
$g\in H(\overline G)$.
\end{lemma}
\begin{proof}
Suppose that $o(g)=\infty$ and that
$g\not\in H(G)$.
Let $H_1$ and $H_2$ be two isomorphic copies of $H$.
Choose a non trivial element $h\in H$
and let $h_1$ and $h_2$ be its copies in~$H_1$ and $H_2$ respectively.
Now define $$\overline G =(G * H_1 * H_2)/  N$$ where
$N$ is the normal subgroup generated by $g^{-1} h_1h_2$.
Then $\overline G\in \K$ by Lemma~\ref{ogog} and moreover $g\in
H(\overline G)$.

If $o(g)=n< \infty$ we first embed $G\leq (G*K) /N$
where $K=\langle x_1, x_2 \,:\, (x_1x_2)^n=1\rangle$
and $N$ is the normal closure of $g^{-1}x_1x_2$.
Then, by the previous lemma $(G*K) /N\in \K$.
Now, since $o(x_1)=o(x_2)=\infty$, we can apply the first case.
\end{proof}

We next show that the class $\K$ is closed under
limits of relatively ``small'' length.
\begin{lemma}
\label{limits-K}
Let $\gamma<\l$ and $\{G_i:  i< \gamma \}$ be an ascending continuous
chain of groups in~$\K$.
Then the union $G=\bigcup\limits_{i<\gamma} G_i$ also belongs to $\K$.
\end{lemma}
\begin{proof}
Let $H'\leq G$ be a subgroup isomorphic to $H$
and let $x\in G$ be such that $[H',x]=1$.
Since $|H'|<\l$ and $\l$ is regular we have $H'\leq G_{i_0}$
for some index $i_0$. We can also suppose that $x\in G_{i_0}$.
Now $G_{i_0}\in \K$, hence $x=1$.
\end{proof}

\begin{lemma}
\label{HG=G}
If $G\in \K$ then there is a group $G'\in \K$ such that
$G \leq G'$ with
 $\card{G'}=\card{G} \cdot \aleph_0$ and $H(G')=G'$.
\end{lemma}
\begin{proof} By Lemma~\ref{free-groups-K} we may assume that $G$ is
infinite, hence the last cardinal condition becomes simply
$\card{G'}=\card{G}$. Applying Lemma~\ref{og} $\card{G}$-times and
Lemma~\ref{limits-K} once we can find an equipotent group $G_1\in
\K$, such that $G\leq G_1$ and  $G \leq H(G_1)$, inductively we find
$G\leq G_{n+1} \in \K$ such that $G_n \leq H(G_{n+1})$ for each
natural $n \geq 1$.\linebreak
The union $G'$ of this chain belongs to $\K$ by
Lemma~\ref{limits-K} and obviously satisfies \linebreak $H(G') = G'$.
\end{proof}

For example, using Lemma~\ref{free-groups-K} and Lemma~\ref{limits-K}
we obtain that, for every cardinal $\k<\l$, the free product
$*_{\a\in \k} \wH_\a$ belongs to $\K$, where $\wH_\a$ is an isomorphic
copy of $\wH$ for every $\a\in \k$.
Now we can apply the following proposition to obtain a group $G'\in \Ks$,
such that $*_{\a\in \k} \wH_\a\leq G'$.

\begin{proposition}
\label{Ks-dense-K}
$\Ks$ is dense in~$\K$.
\end{proposition}
\begin{proof}
Let $G\in \K$. We have to construct a group $G'$ in~$\Ks$ such
that $G\leq G'$ with $|G'|=|G|\cdot \aleph_0$. This can be
obtained easily by using alternatively Lemma~\ref{HNN} and
Lemma~\ref{HG=G} a countable number of times. The limit group will
still lie in~$\K$ by Lemma~\ref{limits-K}. Moreover it is in~$\Ks$
by construction.
\end{proof}

We remark that $\Ks$ is also closed under direct limits
of relatively small length. The proof is similar to that of
Lemma~\ref{limits-K}.

\begin{lemma}
\label{limits-Ks}
Let $\gamma<\l$ and let $\{G_i:  i< \gamma\}$ be an ascending continuous
chain of
groups in~$\Ks$.
Then the union $G=\bigcup\limits_{i<\gamma} G_i$ also belongs to
$\Ks$. \qed
\end{lemma}

\begin{lemma}
\label{amalgamate-Ks}
Let $G^*=G *_{G_0}G'$ be the free product of $G$ and $G'$
amalgamating a common subgroup $G_0$. Suppose that
$G$, $G'$ and $G_0$ are in $\Ks$. Then $G^*\in \K$.
\end{lemma}
\begin{proof}
Let $H'\leq G^*$ be a subgroup isomorphic to $H$, and $x\in G^*$ such
that $[H',x]=1$. As in previous results we can assume that $H'\leq G_0$
and that $x=g_1g_2 \cdots g_n$, written in a reduced form,
has length bigger than two.
Then we have $g_1^*: H' \cong (H')^{g_1}$ both of them inside $G_0$.
Since $G_0\in \Ks$ there exists a $g\in G_0$ such that
$g^*: H'\cong (H')^{g_1}$. We can also suppose that the automorphism
group $\wH'$ is already in $G_0$, by Lemma~\ref{HNN-1}.
Hence, the composition $(g_1^{-1}g)^*: H'\to H'$ is an automorphism,
which is inner by completeness. Thus, $g_1^{-1}g\in G_0$ and thus
$g_1\in G_0$. This is a contradiction, since $x$ was written in a reduced form.
\end{proof}

In order to show that our final group $\G$ is simple we only must
consider normal subgroups $N$ of $\G$ which are cyclically
generated, i.e. there is an $1\neq x\in \G$ with $N=\langle
x^\G\rangle$. We need that $N=\G$. There are two natural cases
depending on the order of $x$. The case that $x$ has infinite
order is taken care by the next Proposition~\ref{Brit}. Hence
assuming that all elements of infinite order are conjugate, a
consequence of Proposition~\ref{Brit}, we only need to note that
any element of finite order can be written as a product of two
elements of infinite order, just take $y$ from a different factor
then $g=(gy) y^{-1}$. Hence $\G=N$. If $x$ has finite order, then
there is a conjugate $y$ of $x$ such that $xy$ has infinite order.
Hence $xy \in N$ and the first case applies. In
Proposition~\ref{densimp} we shall also use this argument.

\begin{proposition}
\label{Brit}
Let $G$ be a group in~$\K$. Let $g, f \in G$, where
$o(f)=o(g)=\infty$ and $g$ does not belong to the normal subgroup
generated by $f$. Then there is a group $\overline{G}\in \K$
such that $G\le \overline{G}$ and $g$ is conjugate to~$f$ in~$\overline{G}$.
\end{proposition}
\begin{proof}
Let $\alpha: \langle f\rangle \to \langle g\rangle$ be the isomorphism
mapping $f$ to $g$. By hypothesis $\alpha \not\in \pinn{G}$.
As in Lemma~\ref{HNN} consider the HNN extension
$\overline{G} =\langle G,t : t^{-1}f t = g \rangle$. We must show
 that $\overline{G} \in \K$. Clearly $\card{\overline{G}} < \l$  and consider
 any $H'$ with $H\cong H'\le G$ and any  $x\in \overline{G}$ with
$[H', x]=1$. By the argument in Lemma~\ref{HNN-2} we may assume that
$H' \leq G$ and $x \in \overline{G}$ with $[H', x]=1$. Now we
apply Lemma \ref{nontrivial-centralizer} as we did in Lemma~\ref{HNN-2}.
\end{proof}



\section{Small cancellation groups}
\label{Section-small-cancellation}

\setcounter{equation}{0}

Let $H$ be as in Section~\ref{a-c}.
Let $\Out(G)= \Aut(G)/\Inn(G)$ be the group of
outer automorphisms of a group $G$. We say that $G$ is {\it complete}
if $\Out(G)=1$ and if the centre $Z(G)=1$, see Robinson \cite{Rob} for
instance. Furthermore, $G$ is called {\it co-hopfian} if every
monomorphism from $G$ to $G$ is an automorphism.

 In this section we shall use small cancellation theory
as in Shelah \cite{She80} or Schupp \cite{Schupp} with the purpose
of ``killing'' monomorphisms which are not inner automorphisms.
This will yield that $\G$ is complete and co-hopfian. Recall from
\cite{Schupp} that an automorphism $\psi$ of a group $G_0$ is
inner if and only if whenever $G_0$ is embedded into a group $G$
then $\psi$ extends to some automorphism of $G$.

An alternative tool would be the use of wreath products, as in
\cite{DG} where Dugas and G\"obel answered a question of Hall
about the existence of particular complete groups.

 Let us first draw a picture which will appear infinitely many
times along the construction of~$\G$:
\begin{equation}
\label{picture}
\begin{array}{rrcll}
 &     & G^*   &       & \\
 & \nearrow &  \downarrow   &\hspace{-0.2cm}\nwarrow    & \\
x_0\in H_0\subset G_1\hspace{-0.2cm} & \rightarrow    &
G &\hspace{-0.2cm} \leftarrow   &G_2
\supset H_1 \ni x_1\\
 & \nwarrow &   & \hspace{-0.2cm}\nearrow & \\
 &    z\in \hspace{-0.3cm} &   G_0    & & \\
\end{array}
\end{equation}
Here, $G_0\leq G_1$, $G_0 \leq G_2$ are all groups in $\K$ and
$G^* =G_1*_{G_0}G_2$.
For $i=0,1$ we have subgroups $H_i\leq G_{i+1}$ isomorphic to $H$
such that $H_i\cap G_0=\{1\}$, and elements $1\neq x_i\in H_i$, $z \in G_0$.

Our task now is to find a new group $G$ in $\K$ where
$G_0\leq  G$ such that:

\medskip
\noindent
(1) $z=\tau(x_0,x_1)$ in~$ G$ for some word
$\tau(x_0,x_1)$ depending on the letters $x_0$ and $x_1$.

\medskip
\noindent (2) $ G$ contains $G_1$ and $G_2$ as malnormal
subgroups, i.e. for every element $y\in G\setminus G_i$ we have
$(G_i)^y \cap G_i = \{1\}$.

\medskip
\noindent
(3)
Suppose that $\psi_0:G_0\to G_0$ is a monomorphism which
extends to a monomorphism $\psi$ with domain $ G$.
Suppose that $\psi \restr H_i=y_i^*\restr H_i$ with
$y_i\in G_{i+1}\setminus C_{G_{i+1}}(x_i)$ for $i=0,1$.
Then $y_0, y_1\in G_0$.

\medskip

If $\psi$ is a monomorphism of our final group $\G$ which is not
an inner automorphism then we will deal with this map restricted
to many subgroups $G$ of $\G$ which are in~$\K^*$. On these
restrictions $\psi\restr G$ we will recognize that $\psi\restr G$
is not inner, by the following result. Call a monomorphism of a
group {\it $\l$-inner\/} if it acts by conjugation on all
subgroups generated by $<\l$ elements. In particular
$\aleph_0$-inner is locally inner in the classical sense.

\begin{proposition}
\label{l-inner}
Suppose that our final group $\G$ is constructed
as a direct limit of groups in $\Ks$ and $\G=H(\G)$.
If $\psi$ is a monomorphism of $\G $ which
is $\l$-inner, then $\psi$ is an inner automorphism.
\end{proposition}
\begin{proof}
Since $\G$ satisfies all conditions for groups in~$ \Ks$ except that
$\card{\G} = \lp$, we know that $H \le \G$ and $\psi(H)$ lies in
some subgroup which is in $\Ks$ hence $\psi\restr H$ is conjugation by
some element in~$\G$. Hence we can suppose without loss of generality
that $\psi\restr H=\ident$.

Now let $H'\leq \G$ be another copy of $H$. As before there is an
element $y\in \G$ such that $\psi \restr H'= y^* \restr H'$. The
subgroup $H'$ and $y$ generate a subgroup $G\leq \G$
which is generated by $<\l$ elements, and it follows from properties of $\G$
that we may also assume $G \in \Ks$.
By hypothesis on $\psi$ there is an element $g\in \G$
such that $\psi \restr H H' =g^*\restr HH'$. Therefore, for all $h\in
H$ and $h'\in H'$, we will have $$ h y^{-1} h' y = \psi(h) \psi(h')=
\psi(hh')= g^{-1} h h' g .$$ In particular, if $h'=1$, then $g^{-1} h
g =h$ for all $h\in H$. Since $G\in \K$ it follows that $g=1$. Thus
$h'=y^{-1} h' y$ for all $h'\in H'$. By the same reason $y=1$. We
have shown that $\psi\restr H' =\ident$ for all $H'< \G$. Since
$H(\G) = \G$, we conclude $\psi= \ident$.  \end{proof}

Here we give an argument
which will provide extensions $G$ of $G_0$ which do not
allow extensions of $\psi\restr G_0$. However $G$ will be
build into $\G$, hence the total map $\psi$ is a contradiction.

More technically, the relations given by (1) will ensure that $z$
can only be moved by conjugation.
This will prove that $\G$ is complete and co-hopfian at the same time.

\medskip
The following ``long'' and ``complicated'' word in $x_0$ and $x_1$
is a good candidate for our purposes:
\begin{equation}
\label{tau}
\tau(x_0,x_1) =  (x_0x_1)(x_0x_1^2)(x_0x_1)^2(x_0x_1^2)^2 \cdots
                 (x_0x_1)^{80} (x_0x_1^2)^{80}.
\end{equation}
Then define
\begin{equation}
\label{r}
r(z,x_0,x_1): = z^{-1} \tau(x_0,x_1).
\end{equation}

Let $N$ be the normal subgroup generated by $r$ in $G^*$. Our
desired $ G$ will be $G^*/N$. Small cancellation theory provides
us a way to decide when a word $w$ in $G^*$ belongs to~$N$. This
is expressed by Proposition~\ref{R-N}. Any weakly cyclically
reduced words in~$N$ must contain relatively ``big parts'' of
words in the symmetrized closure $R$ of $r(z,x_0,x_1)$. We need to
recall the notions involved in this result. More details can be
found in the book by Lyndon and Schupp \cite[p.~286]{LS77} or in
Shelah \cite[p.~378]{She80}.

Let $g=g_1 g_2 \cdots g_m$ be a canonical representation of a word
$g \in G^*=G_1*_{G_0}G_2$. Let $\card{g}=m$ denote the {\it length} of
$g$. Recall that $g$ is {\it weakly cyclically reduced} if $m$ is
even, or $m=1$ or $g_mg_1\not\in G_0$, equivalently, $g_1\cdots g_m$
has no conjugate of length $<m-1$. Note that any word is conjugate
to a weakly cyclically reduced word. Let $R$ be the {\it
symmetrized closure} of $\{r(z,x_0,x_1)\}$ in $G^*$, i.e. add all
inverses and all conjugates which are weakly cyclically reduced.

A {\it piece} (relative to $R$) is an element $b$ in~$G^*$
for which there exist
two different elements $r, r'\in R$ such that $r=bc$ and $r'=bc'$.
Hence, a piece is a subword in~$R$ which cancels out in the multiplication
of two non-inverse elements of $R$.
We say that $R$ satisfies the small cancellation condition $C'(\frac{1}{10})$
if {\it all pieces are relatively small}, that is,
if for some $r\in R$, there is a piece $b$ such that $r=bc$,
then $\card{b} \leq \frac{1}{10} \card{r}$.

It is known from Lyndon and Schupp \cite[p.~289]{LS77} that if $R$
is the symmetrized subset generated by the word $r(z,x_0,x_1)$ of
(\ref{r}), then $R$ satisfies $C'(\frac{1}{10})$. In particular,
Theorem~1.9 in Shelah \cite{She80} yields part 1 of the following
proposition. Part 2 is obvious.

\begin{proposition}
\label{R-N}
Let $G^* = G_1 *_{G_0} G_2$, $R$ and $N$ be as before.
Let $w$ be a weakly cyclically reduced word
of  $G^*\setminus \{1\}$.
Then the following holds:
\begin{enumerate}
\item
If $w\in N$, then
$w$ contains more than $\frac{7}{10}$ of an element $r$ of $R$,
that is, $w$ has a subword $w_1$ which is also contained in~$r$
and such that $\card{w_1} \geq \frac{7}{10} \card{r}$.
\item
If $w\not\in N$, then $w$ does not contain~$\frac{1}{2}$
of any element $r\in R$,
that is, no $r\in R$ has a subword $w_1$ also contained in~$w$
with $\card{w_1} \geq \frac{1}{2} \card{r}$.\qed
\end{enumerate}
\end{proposition}

\begin{proposition}
\label{small-cancellation} Let $G^* = G_1 *_{G_0} G_2$, $R$ and $N$
be as before. Let $ G= G^*/N$. Then $ G\in
\K$, $G_0\leq  G$ and satisfies the properties (1), (2)
and (3) stated before.
\end{proposition}
\begin{proof}
Since elements of $G_1$ and $G_2$ have length 1 in~$G^*$, we have
that $G_1$ and $G_2$ embed in~$ G=G^*/N$, by Proposition~\ref{R-N}. It also
follows that $G_1$ and $G_2$ are malnormal subgroups of
$ G$; see Lyndon and Schupp \cite[p.~203]{LS77} or
Shelah \cite[p.~383]{She80}. Hence $ G$ satisfies
properties (1) and (2).

Let us prove that $ G$ is in $\K$:
Let $H'\leq  G$ be a subgroup isomorphic to $H$ and
$x\in  G$ such that $[H',x]=1$.
We now apply a version of the Torsion Theorem, see
\cite[p.~288, Theorem~11.2 and remarks after the theorem]{LS77}
for a fixed element $h \neq 1$ (of finite order)
in $H'$. Hence, there exists an element $g\in  G$ such that
the conjugate $h^g$ belongs to $G_1$ or $G_2$;
we may assume $h^g \in G_1$. From  $[H', x]=1$ follows
$[h,x] = 1$, hence $[h^g,x^g] = 1$ and
$$1 \neq (h^g)^{x^g} = h^g \in G_1 \cap G_1^{x^g}.$$
Malnormality of $G_1$ in~$ G$ implies that
$x^g$ also belongs to $G_1$. By hypothesis $G_1\in \K$, hence $x^g=1$ and
thus $x=1$.

Finally we show that (3) holds. Let $\psi_0:G_0\to G_0$ be a given monomorphism
which extends to a monomorphism $\psi$ with domain~$ G$.
Suppose that $\psi \restr H_0$ is conjugation
by some element $y_0 \in G_1\setminus G_0$.
We shall assume for simplicity that $y_1=1$
(a similar argument will prove the general case).
Since $\psi$ extends $\psi_0$ we have
$$\tau(\psi(x_0), \psi(x_1))= \psi(z)=\psi_0(z)=g_0$$ for some $g_0 \in G_0$.
In other words, we have the following relation in~$ G$
$$
g_0^{-1} (x_0^{y_0}x_1)(x_0^{y_0}x_1^2)(x_0^{y_0}x_1)^2(x_0^{y_0}x_1^2)^2
\cdots
                 (x_0^{y_0}x_1)^{80} (x_0^{y_0}x_1^2)^{80}=1.
$$ Note that $y_0$ has length 1 in~$G^*$ and that $x_0^{y_0}$
cannot cancel with $x_1$ in~$G^*$. Hence, if we view
$g_0^{-1}y_0^{-1}$ as a new element, then the new word is written
in a weakly cyclically reduced form in~$G^*$.
Proposition~\ref{R-N} tells us that this word must
contain~$\frac{7}{10}$ of some word of $R$. This forces
$x_0^{y_0}=x_0$, that is, $y_0\in C_{G_1}(x_0)$. But this
contradicts the assumption made on $y_0$. Hence $y_0\in G_0$ as
desired.
\end{proof}


\section{A universe for our group constructions.
A uniform partially ordered set of approximations}
\label{upo}

\setcounter{equation}{0}

Let $\lambda$ be as in Main Theorem. We shall need in Lemma~\ref{size-Ku}
below to assume GCH. Let $H$ be a suitable group with automorphism group
$\wH$ of size  $< \l$.

First of all we need to define a ``universe'' where we can embed
our groups. We also want to be sure that there is ``enough room''
for new elements which are added during the construction of our
desired group. This suitable universe is the set
\begin{equation}
\blp := \{\l\a +i : \a\in \lp,\, i<\l\}.
\end{equation}


\noindent
Recall that this is the well-ordered set $\lp \times \l$ of cardinality
$\lp$.
If $x=\l\a+i \in \blp$ with $i<\l$, we write
$\norm{x}=\alpha$ and call $\a$ the {\it norm} of~$x$.
If $X$ is a subset of~$\blp$, the {\it domain} of~$X$,
which we denote by $\dom X$,
is the collection of all~$\norm{x}$ for all $x\in X$;
and its norm $\norm{X}:= \norm{\sup X}$ is the norm
of the supremum $\sup X \in \blp$.
This is well-defined if $\card{X} \leq \l$.

The following picture illustrates how we can embed for instance $H
* H_\a$ in $\blp$, with $\dom(H* H_\a)=\{0, \a\}$. More details
are given in Proposition~\ref{new-Ks-dense-K} below.

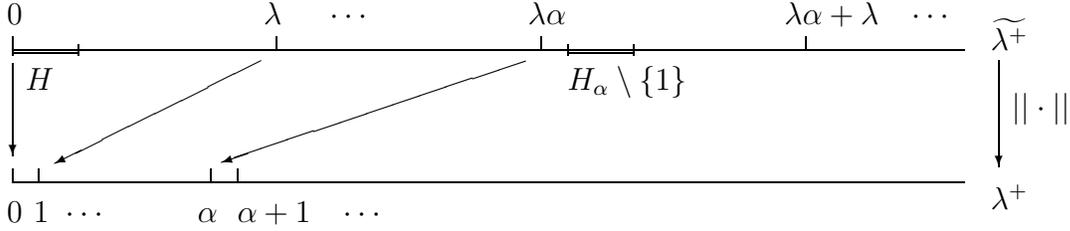
\bigskip

\begin{picture}(300,60)(-10,20)
\put(0,50){\line(1,0){360}}
\put(370,50){\shortstack{$\widetilde{\lambda^+}$}}
\put(373,46){\vector(0,-1){40}}
\put(378,25){\shortstack{$||\cdot ||$}}
\put(0,0){\line(1,0){360}}
\put(370,-10){\shortstack{$\lambda^+$}}

\put(0,45){\vector(0,-1){35}}
\put(94,46){\vector(-2,-1){78}}
\put(194,46){\vector(-3,-1){115}}

\put(0,50){\line(0,1){5}}
\put(100,50){\line(0,1){5}}
\put(200,50){\line(0,1){5}}
\put(300,50){\line(0,1){5}}
\put(0,0){\line(0,1){5}}
\put(10,0){\line(0,1){5}}
\put(75,0){\line(0,1){5}}
\put(85,0){\line(0,1){5}}

\put(95,60){\shortstack{$\l$}}
\put(120,60){\shortstack{$\cdots$}}
\put(195,60){\shortstack{$\l \a$}}
\put(292,60){\shortstack{$\l\a+\l$}}
\put(340,60){\shortstack{$\cdots$}}
\put(-2,60){\shortstack{0}}
\put(-2,-15){\shortstack{$0$}}
\put(8,-15){\shortstack{1}}
\put(20,-15){\shortstack{$\cdots$}}
\put(70,-15){\shortstack{$\a$}}
\put(85,-15){\shortstack{$\a +1$}}
\put(125,-15){\shortstack{$\cdots$}}


\put(210,49){\line(1,0){25}}
\put(210,48){\line(0,1){4}}
\put(235,48){\line(0,1){4}}
\put(210,35){\shortstack{$H_\alpha\setminus\{1\}$}}

\put(0,49){\line(1,0){25}} \put(25,48){\line(0,1){4}}
\put(0,48){\line(0,1){4}} \put(5,35){\shortstack{$H$}}

\end{picture}

\vspace{1.5cm}

\begin{definition}
\label{u-group}
{\rm
Let $u\subset \l^+$ be a subset of cardinality $< \l$ with $0\in u$.
Then a group $G$ of size $|G|<\l$ is called an {\it $u$-group}
if the following holds:
\begin{itemize}
\item[(a)]
$\Dom G \subset \bu := \{ \l\a +i : \a\in u,\, i<\l\}$, where
$\Dom G$ denotes the underlying set of elements of the group~$G$.
We often identify $\Dom G= G$.

\item[(b)]
For every $\d\in \blp$ divisible by $\l$, the subset $G\cap \d$ is
a subgroup of~$G$.
\item[(c)] $\dom G=u$, that is, for every ordinal $\alpha\in u$,
there exists an element $g\in G$
such that $\l\a \leq g < \l(\a +1)$.
\end{itemize}
}
\end{definition}


We now combine Definition~\ref{u-group} with the constructions given in
Section~\ref{a-c}.
Consider two new families of groups $s\K$ and $s\Ks$:

\begin{definition}
{\rm
\begin{enumerate}
\item
A group $G$ belongs to $s\K$ if $G$ is an $u$-group
for some $u\subset\lp$, such that, for every $\d$ divisible by $\l$, the subgroup
$G\cap \delta$ belongs to $\K$ and satisfies the following condition:

If two copies of~$H$ in~$G\cap \delta$ are conjugate by some
element $g\in G$, then $g\in G\cap \delta$.
\item
A group $G$ belongs to $s\Ks$ if $G$ is a member of~$s\K$ such that
for every ordinal $\delta$ divisible by $\lambda$, the
subgroup $G\cap \delta$ belongs to $\Ks$.
\end{enumerate}
}
\end{definition}

Note that $G\in \K$ (resp. $G\in \Ks$) follow respectively
from $G\in s\K$ (resp.  $G\in s\Ks$) and $\card{G}<\l$ because we can
choose $\d\in \blp$ divisible by $\l$ with $G\cap \d=G$.

\begin{notation}
\label{new-K-Ks}
{\rm
We will not need the old families $\K$ and $\Ks$ anymore.
So from now on we rename $s\K$ by $\K$ and $s\Ks$ by $\Ks$, respectively.
}
\end{notation}

For a fixed $u\subset \lp$ with $\card{u}<\l$ we
denote by $\K_u$ (resp.~$\Ks_u$) the set of all $u$-groups in the new
$\K$ (resp.~in~$\Ks$).

\begin{proposition}
\label{new-Ks-dense-K}
$\Ks_u$ is dense in~$\K_u$ for any $u\subset \lp$ with $\card{u}<\l$ and
$0\in u$.
\end{proposition}
\begin{proof}
The proof is the same as that of Proposition~\ref{Ks-dense-K},
but we must indicate how to assign ordinals conveniently to
the constructions described in Section~\ref{a-c}.
In particular we next show that $\K_u$ is not empty.
As $0\in u$, we identify the automorphism group $\wH$ of $H$ with
the first $\card{\wH}$ ordinals in~$\l$; the identity element $1$
goes to $0\in \blp$. Similarly, for every $\alpha\in u$ we
consider an isomorphic copy $\wH_\alpha$ of~$\wH$, and we identify
$\wH_\alpha\setminus \{1\}$ with the first $\card{\wH}-1$ elements
in~$[\l\a,\l\a+\l)$. Let $G_0= *_{\a\in u} \wH_\a$. Hence $G_0$
satisfies part (c) in Definition~\ref{u-group}. If $x\in G_0$ then
$x$ has a unique reduced normal form $x=x_1 x_2 \cdots x_m$, where
each $x_i$ belongs to some factor $\wH_{\a_i}$ and successive
factors are different. Let $\norm{x}$ be the maximum of the norms
$\norm{x_i}$ for $i\leq m$. Now we assign to $x$ any ordinal
in~$[\l \norm{x}, \l\norm{x} +\l)$. There is enough room in this
segment, since $\l$ is uncountable and $|H|\leq |\wH| <\l$, we
simply choose any free place. Hence $G_0$ is an $u$-group. This
yields a new version of Lemmas~\ref{HinK} and \ref{free-groups-K}.
The assignments for Lemma~\ref{limits-K} are clear by
continuity and note that we can do this because the length of the chains
is strictly smaller than $\l$.

Suppose that $G$ and $G'$ have been already embedded in $\subset
\blp$, and let $G^*=G*_{G_0}G'$ for some common subgroup $G_0$. If
$x\in G^*$ consider all possible reduced words $x= w(x_1, ...,
x_n)$ with $x_i\neq 1$ in one of the factors of $G^*$. Assign $x$
to $\blp$ such that $\norm{x}$ is the minimum of all norms of such
words. The same idea works for the rest of results in
Section~\ref{a-c}. In particular, for HNN extensions we assign the
new element $t$ which conjugates $\wH$ and an isomorphic copy
$\wH_1$ to a new ordinal $t\in \bu$ such that
$\norm{t}=\norm{\wH_1}$.
\end{proof}

\begin{lemma}
$\Ks$ is closed under increasing chains of length $<\lambda$.
\end{lemma}
\begin{proof}
The proof is straightforward.
\end{proof}

\begin{definition}
{\rm Let $G_i$ be $u_i$-groups for $i=1,2$. A group isomorphism
$\psi: G_1\to G_2$ is a {\it strong isomorphism} if $\Dom \psi:
\Dom G_1 \to \Dom G_2$ preserves the order in~$\blp$ and it is
constant modulo $\l$, i.e. for every $x\in \Dom G_1$,
$\psi(x)\equiv x\, (mod\, \l)$. This is equivalent to say that if
$x=\l\a +i$ for some $\a\in u_1$ and $i<\l$, then $\psi(x)=\l\a'
+i$ for some $\a'\in u_2$. }
\end{definition}
Clearly, the relation of being strongly isomorphic is an equivalence
relation on
$\K$ and hence on $\Ks$.
Let $[\K]$  denote the set of equivalence classes
$[G]$ for $G\in \K$ and $\pi: \K \to [\K]$ be the natural projection.

\begin{lemma}
\label{size-Ku}
(Assume GCH) For every $u\subset \lp$ with $\card{u}< \l$
and $0\in u$,  the set $[\K_u]$ has cardinality $\leq \l$.
\end{lemma}
\begin{proof}
For a fixed $u\subset \lp$ with  $\card{u}< \l$ we have
exactly $\l^{<\l}$ different subsets in~$\bu$ of size $<\l$,
and $\l^{<\l}=\l$ by GCH.
Now the number of group structures supported by some set
$X$ is $\leq 2^{\mid X\mid}$. By  GCH again and $\card{X}<\l$ we have
$2^{\mid X \mid}=\card{X}^+ \leq \l$. We conclude $\card{[\K_u]}\leq \l$.
\end{proof}

In fact $\card{[\K_u]}=\l$. This follows given other assignments
in the proof of Proposition~\ref{new-Ks-dense-K}. It is
based on a decomposition of~$u$ into $\l$ subsets of size $\l$, which
however will not be needed.

Hence there exists an inclusion (in fact a bijection) $\cod:
[\K_u]\to \l$, which induces a map $\cod: \K_u \stackrel{\pi}{\to}
[\K_u] \to \l$. We can assume that this bijection respects partial
orders, i.e. if $G\leq G'$ then $\cod G\leq \cod G'$.

We finally collect all this information in a {\it coding} map
$$
c: \K\to \l\times \Pa_{<\l}(\l^+), \quad \mbox{given by $G\mapsto (\cod G,
\dom G)$,}
$$
with $c(G)=c(G')$ if $G$ and $G'$ are strongly isomorphic.
Let $\PP : = c(\Ks)$.
Then the partial order on $\K$ induces a partial order on $c(\K)$
and in particular on $\PP$
as follows: we say that $p \leq q$ if there are groups $G_p$ and
$G_q$ in $\K$ coded by $p$ and $q$, respectively, such that $G_p\leq G_q$.
Observe that every element $p=(\alpha,u)$ in~$\PP$
codes every strong isomorphism class $[G]$ in~$\Ks$, hence it is natural to
view the elements of $\PP$ as groups and to consider $\Aut(p)$ for instance.

\begin{theorem} $(\PP, \leq) $ is a $\lp$-uniform partially ordered set.
\end{theorem}
\begin{proof}
We have to check the conditions listed in Definition~A.1.

1. Let $p, q\in \PP$ such that $p\leq q$. Then there are groups
$G_p, G_q\in \Ks$ coded by $p$ and $q$ respectively, such that
$G_p$ is a subgroup of $G_q$. This yields $\dom G_p \subseteq \dom
G_q$, or equivalently $\dom p\subseteq \dom q$, as desired.

2. Let $p, q, r\in \PP$ such that $p,q\leq r$. Let $G_p$,
$G_q$ and $G_r$ groups in~$\Ks$ coded by $p$, $q$ and $r$, respectively,
such that $G_p$, $G_q\leq G_r$.
Consider $G'$ the subgroup of $G_r$ generated by $G_p$ and
$G_q$. It is clear that $\dom G'= \dom p \cup \dom q$
and moreover $G'\in \K$.
Using the fact that $G_r$ is already in~$\Ks$,
we can add to $G'$ the elements of $G_r$,
with norm in the domain of $G'$, in order to obtain a new group
$G''\in \Ks$ such that $G'\leq G''\leq G_r$.
Finally, define $r'= c(G'')$.

The proofs of parts 3 to 6 are straightforward.

7. (Indiscernibility) Suppose that
$p=(\a, u)\in \PP$ and $\varphi: u\to u'$ is an order-isomorphism
in~$\lp$. Let $G_p\in \Ks$ be a group coded by $p$. Using $\varphi$
we can define a strong isomorphism $\psi: G_p \to \psi(G) \subset \blp$ as
follows:
if $x=\l \zeta +i\in G_p$ for some $\zeta\in u\setminus \{0\}$ and $i<\l$,
define $\psi(x)= \l \varphi(\zeta) +i$.
Clearly, the image $\psi(G)$ is a group in~$\Ks$ which has code
$(\alpha, u')=\varphi(p)$. Hence $\varphi(p)\in \PP$.
The fact that $q\leq p$ implies $\varphi(q)\leq \varphi(p)$ is also clear.

8. (Amalgamation property)
%
%
Let $G_p$ and $G_q\in \Ks$ coded by $p$ and $q$ respectively.
Let $G_0\in \Ks$ coded by $p\restr \a$, such that $G_0\leq G_p$
and $G_0\leq G_q$.
Consider $G^*$ the free product of $G_p$ and $G_q$ amalgamating the common
subgroup $G_0$. Then $G^*\in \K$ by Lemma~\ref{amalgamate-Ks}. Now by
Proposition~\ref{new-Ks-dense-K} there is a group $G\in \Ks$ such that
$G^*\leq G$. Take $r=c(G)$.
\end{proof}


\section{Small cancellation groups in~$\blp$}
\label{Section-small-cancellation-blp}

\setcounter{equation}{0}

In this section we formulate a refined version of
Proposition~\ref{small-cancellation}. In the proof of our Main Theorem
we will have to show that the final group $\G$ is complete and co-hopfian.
We will suppose that a monomorphism $\psi:\G\to \G$ exists
which is not an inner automorphism.
Using the help of some prediction principles, we will be able
to find restrictions of $\psi$ on many layers of our group $\G$.
In Proposition~\ref{s-small-cancellation} we provide
a method to get rid of such monomorphism $\psi$.

Recall that $\K$ and $\Ks$ are the new families
introduced in Notation~\ref{new-K-Ks}.

\begin{proposition}
\label{s-small-cancellation} Let $\d_0< \d_1 < \d_2$ be ordinals
in~$\blp$ of cofinality $\lambda$. Suppose that we have $G_i$,
$\psi_i$, $H_i$, $y_i$, $z$ as in the following picture


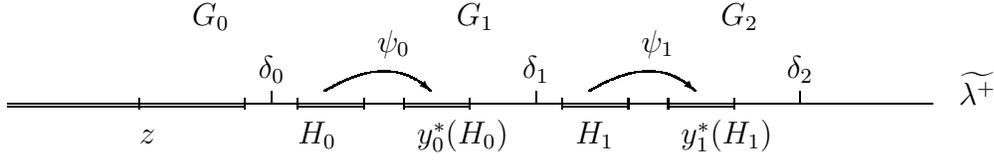
\begin{picture}(300,100)(0,0)

\put(0,50){\line(1,0){350}}
\put(0,49){\line(1,0){90}}
\put(90,48){\line(0,1){4}}
\put(360,50){\shortstack{$\widetilde{\lambda^+}$}}

\put(100,50){\line(0,1){5}}
\put(200,50){\line(0,1){5}}
\put(300,50){\line(0,1){5}}

\put(95,60){\shortstack{$\delta_0$}}
\put(195,60){\shortstack{$\delta_1$}}
\put(295,60){\shortstack{$\delta_2$}}

\put(70,80){\shortstack{$G_0$}}
\put(170,80){\shortstack{$G_1$}}
\put(270,80){\shortstack{$G_2$}}


\put(110,49){\line(1,0){25}}
\put(110,48){\line(0,1){4}}
\put(135,48){\line(0,1){4}}
\put(110,35){\shortstack{$H_0$}}

\put(150,49){\line(1,0){25}}
\put(150,48){\line(0,1){4}}
\put(175,48){\line(0,1){4}}
\put(155,35){\shortstack{$y_0^*(H_0)$}}

\qbezier(120,55)(145,70)(160,55)
\put(160,55){\vector(1,-1){0}}
\put(140,70){\shortstack{$\psi_0$}}

\put(210,49){\line(1,0){25}}
\put(210,48){\line(0,1){4}}
\put(235,48){\line(0,1){4}}
\put(215,35){\shortstack{$H_1$}}

\put(250,49){\line(1,0){25}}
\put(250,48){\line(0,1){4}}
\put(275,48){\line(0,1){4}}
\put(255,35){\shortstack{$y_1^*(H_1)$}}

\qbezier(220,55)(245,70)(260,55)
\put(260,55){\vector(1,-1){0}}
\put(240,70){\shortstack{$\psi_1$}}

\put(50,48){\line(0,1){4}}
\put(50,35){\shortstack{$z$}}

\end{picture}

\vspace{-1cm}

\noindent
where
\begin{enumerate}
\item
$G_i\in \K^*$ for $i=0,1,2$.
\item
$G_i\subseteq \d_i$ and $G_i\cap \d_0=G_0$ for $i=1,2$.
\item
$\psi_i :G_i \to G_i$ are not inner automorphisms such that
$\psi_i\restr G_0= \psi_0$ for $i=1,2$.
\item
$H_0\leq G_1$ and $H_1\leq G_2$ are two copies of~$H$
such that $H_i \cap \d_i= \{1\}$ for $i=0,1$.
\item
$\psi_i\restr H_i= y_i^*\restr H_i$ for some $y_i\in G_{i+1}$ for $i=0,1$.
\item
$z\in G_0$.
\end{enumerate}
Choose $1\neq x_i\in H_i$  such that $\psi_i$
is not conjugation by an element in  $C_{G_{i+1}} (x_i)$ for $i=0,1$
and let $\tau(x_0,x_1)$ be as in (\ref{tau}). Then there is a
small cancellation group $G$ in~$\K$ satisfying the following:
\begin{itemize}
\item[(a)] The relation $z=\tau(x_0, x_1)$ holds in~$G$.
\item[(b)] $G_1$ and $G_2$ are malnormal subgroups of $G$.
\item[(c)]
If $\psi$ is an endomorphism of~$G$  extending $\psi_1$ and
$\psi_2$ and for some $i=0,1$, $\psi \restr H_i$ is conjugation by
some element $y_i\not \in G_0$, then $\tau(\psi(x_0),
\psi(x_1))\not \in G_0$.
\end{itemize}
\end{proposition}

\begin{proof}
Using Proposition~\ref{small-cancellation},
the desired group is the free group $G^*= G_1*_{G_0}G_2$ modulo the
small cancellation $z=\tau(x_0, x_1)$. More precisely $G$ is
the quotient of $G^*$ by the normal subgroup $N$ generated
by the symmetrized closure $R$ of $r=z^{-1}\tau(x_0,x_1)$.
Hence (a), (b) and (c) hold.
\end{proof}


\section{Proof of the Main Theorem}
\label{PROOF}

\setcounter{equation}{0}

In this section we finally prove the existence of the
group $\G$ satisfying the properties listed in the Main Theorem.

Theorem~\ref{HLS-theorem} will ensure that under $\Dl$
player II has a winning strategy for the $\PP$-game,
where $\PP$ is the $\lp$-uniform partially ordered set constructed
in Section~\ref{upo}. In other words player II will be able to
find an increasing and unbounded sequence of ordinals $\zeta_\a$
in~$\lp$ and an ascending chain of admissible ideals $\GG_\a$ on
$\PP_{\zeta_\a}$ which meet the density systems defined by player
I. In the following Subsections 8.1, 8.2 and 8.3 we will tell
player~I how he must define these density systems. These
restrictions will guarantee that our final group $\G$ has
cardinality $\lp$, is complete (d), co-hopfian (e) and simple (f).
Conditions (a), (b) and (c) of the Main Theorem will be
automatically satisfied, since all the groups in the direct system
are in~$\Ks$.

The group $\G$ at the end will be
\begin{equation}
\label{lp-filtration-G}
\G=\bigcup\limits_{\a<\lp} \G_\a,
\end{equation}
where each $\G_\a$ is the union of all groups (coded) in
the directed system $\GG_\a$. Every $\G_\a$ has cardinality $\leq \l$, so that
(\ref{lp-filtration-G}) is a $\lp$-filtration of $\G$.

Let us start playing the $\PP$-game. We proceed by transfinite
induction on $\a< \lp$. By continuity we shall only consider the
induction step from $\a$ to $\a+1$. Suppose that player II has
already found an ordinal $\zeta_\a< \lp$ and an admissible ideal
$\GG_\a$ in~$\PP_{\zeta_\a}$ where $$ \PP_{\zeta_\a}=\{ p \in \PP:
\dom p \subseteq \zeta_\a\}. $$ Now it is the turn of player I. He
chooses $p_{\a+1}\in \PP/\GG_\a=
\{p \in \PP: p\restr\zeta_\a \in \GG \},$ such that
$\norm{p_{\a+1}}> \zeta_\a$. This will guarantee that the length
of the game is $\lp$. He can also define, for each
$\zeta_\a<\b<\lp$, a~list of $(\zeta_\a, \b)$-density systems
$D^i$ over $\GG_\a$, where $i$ ranges over a set of cardinality
$\leq \l$.

Given $u\subseteq v \subseteq \lp$ with $\card{v}< \l$, we denote
for the rest of this section $$\E:=\{ p\in \PP / \GG_\a : \mbox{
$\dom p \subseteq v \cup \b$}\}.$$ Next we define $\l$ subsets
$D(u,v)$ in~$\E$ satisfying the properties stated in
Definition~\ref{density-system}. In particular, property~2 follows
directly if we define in the same way $D(u',v')$ for all pairs
$(u',v')$ such that $u \cap \beta = u' \cap \beta$, $v \cap \beta
= v' \cap \beta$ and there is an order-isomorphism from $v$ onto
$v'$ which maps $u$ onto $u'$.


\subsection{Density system to obtain~$\card{\G}=\lp$}


We want that our group $\G$ has $\dom \G=\l^+$.
This will imply that $\card{\G}=\lp$.
With this objective in mind player I chooses in his $\a th$ move
the following density system:

If $u, v $ are subsets in~$\lp$ such that $u\subseteq v$ and $\card{v}<\l$
player I defines
$$
D(u,v) :=\{ p\in \E : \mbox{ if $\a \not \in u$ then $\a\in \dom  p$}\}.
$$

\begin{proposition} \label{dendom}
$D$ is a density system over $\GG_\a$.
\end{proposition}
\begin{proof}
It is clear that $D(u,v)$ is closed upwards in $\E$.
Now we have to show that $D(u,v)$ is dense in~$\E$.
In the case when $\a\in u$ then $D(u,v)= \E$ and density is obvious.
Suppose that $\a\not\in u$ and let $q$ be any element in~$\E$.
Let $G$ be an $u$-group which is coded by $q$.
As in the proof of Proposition~\ref{new-Ks-dense-K} we can find a group
$G'\in \Ks$, such that $G\leq G'$,
which contains a subgroup $H'$ isomorphic to $H$ and
$\a \in \dom H'\subseteq [\l \a, \l \a +\l)$.
In particular, we have $\a\in \dom G'$.
Moreover, the group $G'$ can be constructed in such a way that
$\dom G'\subseteq v \cup \b$.
Therefore, $p:=(\cod G', \dom G')\in D(u,v)$ and $q\leq p$ as desired.
\end{proof}



\subsection{Density systems to obtain that $\G$ is complete and co-hopfian}


Suppose for a moment that $\G$ is already constructed as a limit of
groups in $\Ks$ and $\psi$ is an unwanted monomorphism
of $\G$ which is not a inner automorphism. Hence, it is not inner
on many layers by Proposition~\ref{l-inner}.

First we note that for a fixed pair $(u,v)$ and $\zeta_\a<\beta< \lp$
the set $E$ of equivalence classes of all monomorphisms `up to order
isomorphism' of $u$-groups $U$ with $u \subseteq v \cup \beta$ has
size~$\l$. Choose then an enumeration of $U \times E$ by $\l$. For
every $i\in\l$ we shall define a $(\zeta_\alpha,\beta)$-density
system $D^i$ taking care of an element $z \in U$ and the monomorphism
$\psi$ coded by $i \in \l$. If $p \in \E$ then consider $p \restr u = G_0$
and let $z \in G_0, \psi \in \Aut(G_0)$ be the pair of elements coded
by $i < \l$. Then we put $p \in D^i (u,v)$ if $\psi$ is not inner and
$G_0$ extends to $p = G$ with $G_1, G_2, \psi_1, \psi_2, H_0$ and
$H_1$ as in Proposition~\ref{s-small-cancellation}.

Note that such a monomorphism of $G_0$ can not survive to a
monomorphism of $G$. On the other hand if $p\in \E$ such that $\psi$
is inner on $p\restr u$ then $p\not\in D^i(u,v)$.

\pagebreak

\begin{proposition}\label{dencomp}
$D^i $ is a density system over $\GG_\a$ for every $i\in \l$.
\end{proposition}
\begin{proof}
It is clear that every $D^i(u,v)$ is closed upwards. The fact
that $D^i(u,v)$ is dense in~$\E$ is shown in
Proposition~\ref{s-small-cancellation}, where we can find
a small cancellation group $G\in \K$ with the required properties.
If necessary, use also Proposition~\ref{new-Ks-dense-K} to obtain
a group $G_3$ in $\Ks$ such that $G\leq G_3$.
\end{proof}

\subsection{Density systems to obtain that $\G$ is simple}

Finally we take care of simplicity of our group $\G$. As in the
previous subsection, we shall define $\l$ density systems $D^i$. As
before, let $u\subseteq v \subseteq \lp$ such that $\card{v}<\l$.
Suppose that $U$ is a group in~$\Ks$ with $\dom U=u$. We consider the
set $P_U$ of all ordered pairs $(x,y)$ where $1 \neq y, x \in U$ and
choose a coding map $ c : P_U \hookrightarrow \l.$ Then, for every
$i\in c (P_U)$ player I chooses $D^i (u,v)$ to be $$ \{ p\in \E :
p= c (G), \; p\restr u = c(U), \; c(x,y)=i, \;  {\rm and} \;
x\not\in \langle y^U \rangle \; \Longrightarrow \; x\in \langle y^G
\rangle \}, $$ where $\langle y^U \rangle = \langle y^z :  z\in
U\rangle$.

\begin{proposition} \label{densimp}
$D^i$ is a density system over $\GG_\a$ for every $i$.
\end{proposition}
\begin{proof}
Let $q \in \E$ represent the group $U$ and let $1 \neq y,x \in U$ be
such that  $c(x,y) = i$. If $x  \in \langle y^U\rangle $, then
$D^i(u,v) = \E$ and density is obvious. Suppose that
$x \not \in \langle y^U \rangle $. We have to construct a $v$-group $G$
represented by $p \in \E$ extending $q$ such that $x \in \langle y^G \rangle$.

We consider the three possible cases:

Case 1: $o(x)=o(y)=\infty$. In this case we use
Proposition~\ref{Brit} as follows. Define $G_1 = \langle U, t :
x=t^{-1}y t\rangle$, where $t$ is a new element in $v$ with
$||t||=||x||$. Of course, $G_1\in \K$, so we can construct $G_1
\leq G$ with $G\in \K^*$ and $\dom(G)\subset v$, as explained in
the proof of Proposition~\ref{new-Ks-dense-K}. Hence $x$ and $y$
are conjugate in $G$.

Case 2: $o(x)=n< \infty$ and $o(y)=\infty$. Here, we use
Lemma~\ref{og} (part 2 of its proof). We can find an element $x' =
xz \in G\in \K^*$ of infinite order such that $z^{-1}y$ has
infinite order as well. From case 1 we note that $y,\  xz \mbox{
and } z^{-1}$ are conjugate in $G$. Hence $x \in \langle y^G
\rangle$.

Case 3: If $y$ has finite order, we argue similarly. Take a
conjugate $z$ of $y$ in some $G\in \K^*$ such that $yz$ has
infinite order. Hence $\langle y^G \rangle$ contains all elements
of infinite order and case 1 and 2 apply.

We conclude in any case that $D^i$ is a density system.
\end{proof}

\subsection{Existence of $\G$}
\label{existence-G}
Using the density systems and the winning strategy of the game
Theorem \ref{HLS-theorem} we find the building blocks for the desired
group.

\begin{theorem}
\label{direct-system} (Assume $\Dl$) There is a chain of admissible
ideals $\GG_a \subseteq \PP_{\zeta_\a}$ for an unbounded sequence
$\{\zeta_\a < \lp: \a < \lp\}$ such that $\GG =
\bigcup\limits_{\alpha < \lp} \GG_\a$ satisfies the following
conditions.
\begin{enumerate}
\item
For every $\a\in \lp$ there is a group $G\in \GG$ such that
$\alpha\in \dom G$.

\item Completeness and co-hopfian:
If  $\delta_0<\blp$ divisible by $\l$, $z \in G_0 \in \GG$ with $G_0\subseteq \delta_0$
and $\psi_0:G_0\to G_0$ is monomorphism which is not inner, then there
exist ordinals $\d_0 < \d_1 < \d_2$ in~$\blp$ divisible by $\l$  and groups $G_0 \leq
G_i \in \GG $ with $G_i \subseteq \d_i$ and subgroups $H_0 \leq G_1$,
$H_1\leq G_2$ isomorphic to $H$ such that $H_i\cap
\d_i=  \{1\}$ for $i=0,1$  as in the hypothesis
of Proposition~\ref{s-small-cancellation}.
\item
Simplicity: For every $1 \neq y, x \in G \in \GG$ there is $G' \in
\GG$ such that $G \leq G' $ and  $x \in \langle y^{G'}\rangle$.
\end{enumerate}
\end{theorem}
\begin{proof} We apply Theorem \ref{HLS-theorem} and let player I choose the
density systems above. Condition 1 of the Theorem follows
immediately from Proposition~\ref{dendom}, condition 2 follows
from Proposition~\ref{dencomp} and condition 3 follows from
Proposition~\ref{densimp}.
\end{proof}

Let $\G$ be the group which is the union of the directed system $\GG$
given in Theorem~\ref{direct-system}.
We finally check that $\G$ has the desired properties:

\noindent
{\bf
$\G$ has cardinality $\lambda^+$:}
By Theorem~\ref{direct-system} we have that for every $\a<\lp$, there
is a group $G$ in $\GG$ such that $\a\in \dom G$. This says that
$\dom \G= \lp$ and hence, $\lp \leq \card{\G}$. On the other hand,
all groups $G\in \GG$ lie in  $\blp$. Thus $\card{\G}\leq
\card{\blp}=\lp$ and $\card{\G}=\lp $ follows.

\noindent
{\bf
Properties (a)-(c):}
These properties are satisfied by all the groups $G$
in~$\GG$ by construction. Now if $x\in \G$ and $H_1\leq \G$ is
isomorphic to $H$ such that $[H_1,x]=1$ there exists a group $G\in
\GG$ containing $x$ and $H_1$. Hence $x=1$ from $G \in \Ks$. This
shows property (c). The properties (a) and (b) can be proved
similarly.

\noindent {\bf Properties (d)-(e):} Let $\psi: \G\to \G$ be a
monomorphism (possibly not inner automorphism). Let $C$ be the set
of ordinals $\delta<\blp$ which are divisible by $\lambda$ and
such that $\psi$ maps $\delta \cap \G$ into $\delta \cap \G$. Then
$C$ is clearly a club of~$\blp$. Now let $S=\{\d\in \blp : \cf
(\d)=\l\}$. Then $S$ is stationary in~$\blp$ and hence $S \cap C
\neq \emptyset$.

Choose for each $\d\in S $ a subgroup $H^0_\d$ of $\G$ isomorphic
to~$H$ such that $$H^0_\d \cap \d =  \{1\}.$$
By property (b) we find $y_\d \in \G$ such that
$\psi \restr H^0_\d = y_\d^* \restr H^0_\d$.

Now we have two cases:

\noindent
{\it Case 1:} Suppose that $y_\d < \d$ for all $\d\in S$. In that case we
have a regressive function $$ S \to \blp, \quad \mbox{given by $\d
\mapsto y_\d$}. $$ Hence, by Fodor's lemma (see Lemma~\ref{Fodor}),
there exists a stationary subset $S'\subset S$ and $y\in \G$ such
that $y_\delta=y$ for all $\delta\in S'$.

Let $z\in \G$. We shall prove that $\psi(z)=y^*(z)$. For every
$\alpha\in \lambda^+$ choose an ordinal $\beta=\beta(\alpha)>\alpha$
such that $\l\b\in S'$. Such a $\beta$ exists since $S'$ is
stationary and $\{\b\in \lp: \b\geq \alpha\}$ is a club. Let
$G_\alpha\in \GG$ such that $z\in G_\a$ and $H_\alpha= H^0_{\beta}
\leq  G_\alpha$. Let $H'_\alpha=\psi(H_\alpha)$. Note that
$H'_\alpha$ is a subgroup of~$G_\alpha$ such that $(\dom
H'_\alpha)\cap \a=\emptyset$. Now by Theorem~\ref{direct-system} we
can find ordinals $\alpha_0<\alpha_1<\alpha_2$, groups $G_0$, $G_1$
and $G_2$ in~$\GG$ satisfying the hypothesis of
Proposition~\ref{s-small-cancellation} (take $\d_i=\l\a_i$). Then
there exists a group $G\in \GG$ extending $G_1$ and $G_2$ where the
small cancellation $z=\tau(x_0,x_1)$ holds, where $z\in G_0$ and
$x_i\in H_i$ for $i=0,1$. Hence $$\psi(z)= \psi( \tau(x_0,x_1)) =
\tau (\psi(x_0) , \psi (x_1))= \tau (y^*(x_0), y^*(x_1))=
y^*\tau(x_0,x_1)=y^*(z).$$ Since $z$ is arbitrary, we have shown that
$\psi$ is a inner automorphism.

\noindent
{\it Case 2:} Suppose that for some $\delta\in S$, we have $\psi\restr
H^0_\alpha \neq y^*\restr H^0_\alpha$ for any $y\in \delta$.
We argue as in the case~1 with $z=1$. Hence, there exist ordinals
$\alpha_0<\alpha_1<\alpha_2$ and groups $G_i\in \Ks$ as in
Proposition~\ref{s-small-cancellation} so that $1=z=\tau(x_0,x_1)$
for some elements $x_0$ and $x_1$. Hence $1= \psi(z)=\tau (\psi(x_0),
\psi(x_1)) \not \in G_0$, which is impossible.

Thus $\G$ is complete and co-hopfian, as desired.

\noindent
{\bf
Property (f):
}
This follows from the third condition of Theorem~\ref{direct-system}.

\appendix

\section{The combinatorial principle $\Dl$}
\label{appendix}
\setcounter{equation}{0}

In this appendix we recall for convenience of the reader the
combinatorial result proved by Hart, Laflamme and Shelah in
\cite{HLS93}. This result will guarantee in
Theorem~\ref{direct-system} the existence of a certain direct system
$\GG$ and $\G$ in Main Theorem will be the limit of the groups coded
by elements in $\GG$.

Theorem~\ref{HLS-theorem} is stated in terms of a game of~$\lp$ moves
between two players I and II. In the turn $\a<\lp$ player II picks an
ordinal $\zeta_\a<\lp$ and an admissible ideal $\GG_\a$ from
a partial uniform set $\PP$ (see definitions below).

The elements of~$\PP$ are pairs $p=(\a, u)\in \l \times
\Pa_{<\lambda}(\lp)$. We write $u=\dom p$ and call $u$ the {\it
domain} of~$p$. The ordinals $\a$ are there to code
all possible algebraic structures that $u$ can support.
In our case $\a$ will run over all
strongly isomorphic classes of groups supported by a fixed~$u$ (see
Section~\ref{upo}).

\begin{definition}
{\rm A {\it standard $\lp$-uniform partial order} is a partial order
$\leq$ defined on a subset $\PP$ of~$\l \times \Pa_{<\lambda}(\lp)$
satisfying the following conditions:

\begin{enumerate}
\item {\rm (Compatibility of the orders)} If $ p \leq q $ then $\dom p
\subseteq
\dom q$.
\item For all $p,q,r \in \PP$ with $p,q \leq r$ there is $r' \in \PP$ such that
$p,q \leq r' \leq r$ and $\dom r' = \dom p \cup \dom q.$
\item If $\{p_\a: \a < \delta\}$ is an increasing sequence
in~$\PP$ of length $\delta<\lambda$ then it has a least upper bound
$q \in \PP$ with $\dom q = \bigcup\limits_{\alpha < \delta} \dom
p_\alpha$; we say that $q = \bigcup\limits_{\alpha < \delta}
p_\alpha$.
\item  If $p \in \PP $ and $\alpha < \lambda^+$ then there  is $q \in \PP$
such that $ q \leq p$ and $\dom q = \dom p \cap \alpha $ and there is
a unique maximal such $q$ for which we write $q = p\restr \alpha$.
\item {\rm (Continuity I)} If $\delta$ is a limit then $p\restr \delta =
\bigcup\limits_{\alpha < \delta} p\restr \alpha$.
\item {\rm (Continuity II)}
If $\{p_\a: \a < \delta\}$ is an increasing sequence in~$\PP$ of
length $\d<\a$ then $$ (\bigcup\limits_{i < \delta} p_i)\restr \alpha
= \bigcup\limits_{i < \delta} (p_i\restr \alpha).$$
\item {\rm (Indiscernibility)} If $p = (\alpha, u) \in \PP$ and $\varphi: u
\to u' $
is an order-isomorphism in~$\lp$, then $\varphi(p) := (\alpha,
\varphi(u)) \in \PP$. Moreover, if $q \leq p$ then $\varphi(q) \leq
\varphi (p)$.
\item {\rm(Amalgamation property)} For every $p,q \in \PP$ and $\alpha <
\lambda^+$ with
$p\restr\alpha \leq q$ and $\dom p \cap \dom q = \dom p \cap \alpha$
there is an $r \in \PP $ such that $p,q \leq r$.
\end{enumerate}
If all these conditions hold, $\PP = (\PP,\leq)$ is called a {\it
$\lp$-uniform partially ordered set}. }
\end{definition}

Consider the following filtration of $\PP$, where $\alpha <
\lambda^+$, $$\PP_\alpha := \{ p \in \PP: \dom p \subseteq \alpha
\}.$$

\begin{definition}
{\rm Let $(\PP, \leq)$ be a $\lp$-uniform partially ordered set and
$\a<\lp$. A subset $\GG \subseteq \PP_\alpha$ is an {\em admissible
ideal of~$\PP_\alpha$} if the following holds:
\begin{enumerate}
\item $\GG$ is closed downward, i.e. if $p\in \GG$ and $q\leq p$ then
$q\in \GG$.
\item $\GG$ is $\l$-directed, i.e.
if ${\cal A} \subset \GG$ and $\card{{\cal A}}<\l$ then ${\cal A}$
has an upper bound in~$\GG$.
\item
(Maximality) If $p\in \PP$ is compatible with all $q\in \GG$ then
$p\in \GG$. (Recall that $p$ is compatible with $q$ if there is an
$r\in \PP$ such that $p$, $q\leq r$.)
\end{enumerate}
For an admissible ideal $\GG$ of~$\PP_\alpha$, we define
$$\PP/\GG :=\{p \in \PP: p\restr\alpha \in \GG \}.$$ }
\end{definition}
This consists of all extensions of elements in~$\GG$. Note that
this notion is compatible with taking direct limits in the
following sense. Let $\{\zeta_\beta : \b<\a\}$ be an increasing
sequence of ordinals in~$\lp$ converging to~$\zeta$. Suppose that
for every $\beta$ we have an admissible ideal $\GG_\beta$
of~$\PP_{\zeta_\beta}$. Then there is a unique minimal admissible
ideal containing the set theoretic union of the $\GG_\beta$'s.
With a slight misuse of notation we write $$\GG_{<\a}:=
\bigcup\limits_{\b<\a} \GG_\beta,$$ for this admissible ideal
of~$\PP_{\zeta}$, see Hart, Laflamme and Shelah \cite[p.~173,
Lemma~1.3]{HLS93}.

\begin{definition}
\label{density-system} {\rm Let $\GG$ be an admissible ideal
of~$\PP_\alpha$ and $\alpha < \beta < \lambda^+$. An {\em
$(\alpha,\beta)$-density system over $\GG$} is a function $$D: \{
(u,v) : u\subseteq v\in \Pa_{<\l}(\lp)\} \arr \Pa (\PP)$$ such that
the following holds:
\begin{enumerate}
\item $D(u,v) \subseteq \{p \in \PP/\GG: \dom p
\subseteq v \cup \beta \}$ is a dense and upward-closed subset.
\item If  $(u,v),(u',v')$ and $u \cap \beta = u' \cap \beta,\
v \cap \beta = v' \cap \beta$ and there is an order-isomorphism from
$v$ onto $v'$ which maps $u$ onto $u'$, then for any ordinal $\gamma$
we have $$(\gamma,v) \in D(u,v) \iff (\gamma,v') \in D(u',v').$$
\end{enumerate}

An admissible ideal $\GG'$ of~$\PP_{\alpha'}$ for some $\alpha' <
\lambda^+$ {\em meets} the $(\alpha,\beta)$-density system $D$ (over
$\GG$) if $\alpha < \alpha', \ \GG \subseteq \GG'$ and for each $u
\in \Pa_{<\l}(\a')$ there is a $v \in \Pa_{<\l}(\a')$, with $u
\subseteq v$ and such that $D(u,v) \cap \GG' \neq \emptyset$. }
\end{definition}

\bigskip
\noindent {\bf The $\PP$-game}

\begin{enumerate}
\item
Player I moves first.

\noindent Let $\a$ be a successor ordinal and suppose that player I
has played his $(\a-1)th$ move.

\item Then player II picks an ordinal $\zeta_{\a} < \lp$ and an
admissible ideal $\GG_\a$ of $\PP_{\zeta_\a}$.

\item Then player I chooses an element $p_{\a+1}\in \PP/\GG_{\a}$
and chooses, for every $\zeta_\a < \b <\lp$, a list of $(\zeta_\a,
\b)$-density systems $D^i$ over $\GG_{\a}$, where $i\in I_\a$ with
$|I_\a|\leq \l$.

\noindent If $\a=\lim_{\a'<\a} \a'$ is a limit ordinal the moves go
as follow:

\item   At the $\alpha th$ move, player II will have chosen an increasing
sequence of ordinals $\zeta_{\a'} < \lambda^+ $ and
 will have defined an admissible ideal $\GG_{\a'}$ of~$\PP_{\zeta_{\a'}}$,
for all $\a'<\a$.

\item
   Then player I will choose an element $p_\a\in \PP/\GG_{<\a}$
and he will also choose at most $\l$ density systems $D^i$ over
$\GG_{<\alpha}$, where $i \in I_\alpha$ with~$\card{I_\a}\leq \l$.
\item
   After player I has finished his $\a th$ move player II will
pick again an ordinal $\zeta_\alpha$ and an admissible ideal
$\GG_\alpha$ of~$\PP_{\zeta_\alpha}$.
\end{enumerate}

Player II {\it has a winning strategy} for the $\PP$-game if the
sequences $\zeta_\alpha$ and $\GG_\alpha$ are increasing, we have
$p_\alpha \in \GG_\alpha$ and for all $\b \geq \a$, $\GG_\beta$ {\it
meets} $D_\a^i$ for all $i \in I_\alpha$.

\medskip
\noindent Theorem~\ref{HLS-theorem} below states that player II has a
winning strategy if $\Dl$ holds. This is a combinatorial principle
which helps player~II to predict the moves of its opponent player.

\medskip
\noindent {\it $\Dl$ \, asserts that there are sets ${\cal A}_\alpha
\subseteq \Pa(\alpha)$ for all $\a< \l$, with $|{\cal A}_\alpha|<
\lambda$ such that for all $X \subseteq \lambda$ the set $\{ \alpha
\in \lambda: X \cap \alpha \in {\cal A}_\alpha \}$ is stationary
in~$\lambda$. }

\begin{theorem}
\label{HLS-theorem} If $\Dl$ holds then player II has a winning
strategy for the $\PP$-game.
\end{theorem}

We finally remark that $\lambda = \aleph_0$ or $\Diam (\l)$
implies $\Dl$. Recall that $\Diamond(\l)$ is the particular case
of $\Dl$ where each ${\cal A}_\alpha$ consists on one single
subset of $\l$. Moreover, if $\l=\kappa^+$ with $\cf(\kappa)>
\aleph_0$ and GCH holds, then $\Dl$ holds. We also note that $\Dl$
implies $\l^{<\l}=\l$.


\vskip 0.5 cm

\setlength{\baselineskip}{0.6cm}

\bigskip\noindent
Fachbereich 6, Mathematik und Informatik, Universit\"at Essen,\\
D-45117 Essen, Germany,  e-mail: {\tt R.Goebel@Uni-Essen.De}

\bigskip\noindent
Departament de Matem\`atiques, Universitat Aut\`onoma de
Barcelona,\newline E--08193 Bellaterra, Spain, e-mail: {\tt
jlrodri@mat.uab.es}

\bigskip\noindent
Institute of Mathematics, Hebrew University,\\ 91904 Jerusalem,
Israel, e-mail: {\tt shelah@math.huji.ac.il}

\end{document}